\documentclass[UTF-8,reqno]{amsart}
\usepackage{enumerate, bbm}
\usepackage{amssymb,url,color, booktabs}
\usepackage{tikz}
\usepackage{mathrsfs}
\usepackage{dutchcal}
\usepackage{threeparttable}
\usepackage{color}
\usepackage[colorlinks=true]{hyperref}
\hypersetup{
    linkcolor=blue,          
    citecolor=red,        
    filecolor=blue,      
    urlcolor=cyan
}

\usepackage{color}
\usepackage{ulem}

\definecolor{MyDarkBlue}{cmyk}{0.8,0.3,0.8,0.4}
\definecolor{yellow}{rgb}{0.99,0.99,0.70}
\definecolor{white}{rgb}{1.0,1.0,1.0}
\definecolor{black}{rgb}{0.00,0.00,0.00}

\newcommand{\blue}{\color{blue}}

\numberwithin{equation}{section}

\newcommand{\be}{\begin{eqnarray}}
\newcommand{\ee}{\end{eqnarray}}
\newcommand{\ce}{\begin{eqnarray*}}
\newcommand{\de}{\end{eqnarray*}}
\newtheorem{theorem}{Theorem}[section]
\newtheorem{lemma}[theorem]{Lemma}
\newtheorem{remark}[theorem]{Remark}
\newtheorem{definition}[theorem]{Definition}
\newtheorem{proposition}[theorem]{Proposition}

\newtheorem{corollary}[theorem]{Corollary}

\usepackage[nobysame,abbrev]{amsrefs}
\BibSpec{article}{%
+{}{\PrintAuthors} {author}
+{,}{ \textrm} {title}
+{.}{ \textit} {journal}
+{,}{ \textbf} {volume}
+{}{ \parenthesize} {date}
+{,}{ } {pages}
+{.}{ arXiv:} {eprint}
+{.}{} {transition}
}

\def\bt{\begin{theorem}}
\def\et{\end{theorem}}
\def\bl{\begin{lemma}}
\def\el{\end{lemma}}
\def\bd{\begin{definition}}
\def\ed{\end{definition}}
\def\bp{\begin{proposition}}
\def\ep{\end{proposition}}
\def\bc{\begin{corollary}}
\def\ec{\end{corollary}}
\def\br{\begin{remark}}
\def\er{\end{remark}}
\def\bexa{\begin{example}}
\def\eexa{\end{example}}

\def\bC{{\mathbf C}}

\def\bE{{\mathbf E}}

\def\bP{{\mathbf P}}

\def\bX{{\boldsymbol{X}}}

\def\1{{\mathbf{1}}}

\def\cB{{\mathcal B}}

\def\cS{{\mathcal S}}

\def\mC{{\mathbb C}}
\def\mD{{\mathbb D}}
\def\mE{{\mathbb E}}

\def\mI{{\mathbb I}}

\def\mK{{\mathbb K}}
\def\mL{{\mathbb L}}

\def\mN{{\mathbb N}}

\def\mP{{\mathbb P}}
\def\mQ{{\mathbb Q}}
\def\mR{{\mathbb R}}

\def\mT{{\mathbb T}}

\def\sF{{\mathscr F}}

\def\sI{{\mathscr I}}

\def\sL{{\mathscr L}}

\def\eps{\varepsilon}
\def\t{\tau}
\def\p{\partial}
\def\de{\delta}
\def\g{\gamma}
\def\l{\lambda}
\def\k{\kappa}

\def\[{{\Big[}}
\def\]{{\Big]}}
\def\<{{\langle}}
\def\>{{\rangle}}
\def\({{\big(}}
\def\){{\big)}}

\def\tr{\mathrm {tr}}

\def\dif{{\mathord{{\rm d}}}}

\def\bbX{{\boldsymbol{X}}}

\def\bbx{{\boldsymbol{x}}}
\def\bby{{\boldsymbol{y}}}
\def\bbz{{\boldsymbol{z}}}
\def\bbp{{\boldsymbol{p}}}

\def\bb2{{\boldsymbol{2}}}
\def\no{\nonumber}
\def\={&\!\!=\!\!&}

\def\geq{\geqslant}
\def\leq{\leqslant}
\def\ge{\geqslant}
\def\le{\leqslant}

\def\S{\Sigma}
\def\c{\mathord{{\bf c}}}

\def\eps{\varepsilon}
\def\t{\tau}
\def\e{\mathrm{e}}

\def\t{\theta}

\def\p{\partial}

\def\g{\gamma}
\def\l{\lambda}

\def\[{{\Big[}}
\def\]{{\Big]}}
\def\<{{\langle}}
\def\>{{\rangle}}

\def\tr{\mathrm {tr}}

\def\c{\mathord{{\bf c}}}

\def\bX{{\mathbf X}}

 \def\R{\mathbb R}
 \def\R{\mathbb R}    
\def\N{\mathbb N}  
   
\def\<{\langle} \def\>{\rangle}

\allowdisplaybreaks
\begin{document}
\title[Heat kernel estimates for kinetic SDEs]
{Heat kernel estimates for kinetic SDEs with  drifts being unbounded and  in Kato's class}

\author{Chongyang Ren and Xicheng Zhang}



\address{Chongyang Ren: 
School of Mathematics and Statistics, Wuhan University, Wuhan, Hubei 430072, China\\
Email: rcy.math@whu.edu.cn}

\address{Xicheng Zhang:
School of Mathematics and Statistics, Beijing Institute of Technology, Beijing 100081, China\\
Email: XichengZhang@gmail.com
 }

\thanks{
This work is partially supported by NNSFC grants of China (Nos. 12131019), and the German Research Foundation (DFG) through the 
Collaborative Research Centre(CRC) 1283/2 2021 - 317210226 ``Taming uncertainty and profiting from randomness and low regularity in analysis, stochastics and their applications".}

\begin{abstract}
In this paper we investigate the existence and uniqueness of weak solutions for kinetic stochastic differential equations with H\"older diffusion and unbounded singular drifts in Kato's class.  
Moreover, we also establish sharp two-sided estimates for the density of the solution. In particular, the drift $b$ can be in the mixed $L^q_tL^{p_1}_{x_1}L^{p_2}_{x_2}$ space with
$\frac2q+\frac{d}{p_1}+\frac{3d}{p_2}<1$. As an application, we show the existence and uniqueness of weak solution to
a second order singular interacting particle system in ${\mathbb R}^{d N}$.
\end{abstract}
\maketitle

\section{Introduction}
In this paper we aim to prove the Aronson-like bounds for the transition probability density of the kinetic stochastic differential equations
(abbreviated as kinetic SDEs) in $\mR^{2d}$ with H\"older diffusion and unbounded  singular drifts, 
which are in certain Kato's class:
\begin{align}\label{1}
\left\{
\begin{aligned}
&\dif X^{1}_t=b\(t,X^{1}_t,X^{2}_t\)\dif t+\sqrt{2}\sigma\(t,X^{1}_t,X^{2}_t\)\dif W_t,\\
&\dif X^{2}_t=X^{1}_t\dif t,
\end{aligned}
\right.
\end{align}
where $(W_t)_{t\ge 0}$ is a $d$-dimensional Brownian motion on some stochastic basis  $\mathfrak{F}:=(\Omega,\sF, (\sF_t)_{t\geq 0},\bP)$,
and $(b,\sigma):\mR_+\times\mR^{2d}\to(\mR^d, \mR^d\otimes\mR^d)$ are Borel measurable functions.
For simplicity of notations, we shall denote $\bbX:=(X^1, X^2)$ and for $\bbx=(x_1,x_2)\in\mR^{2d}$,
\begin{align}\label{S}
B(t,\bbx):=(b(t,\bbx),x_1)^*,\ \ \Sigma(t,\bbx) :=(\sqrt{2}\sigma(t,\bbx),0_{d\times d})^*,
\end{align}
where $*$ stands for the transpose of a row vector.
By this, we can rewrite SDE $\eqref{1}$ as
\begin{align}\label{25}
\dif\bbX_t=B(t,\bbX_t)\dif t+\Sigma(t,\bbX_t)\dif W_t,
\end{align}
and its generator is given by
\begin{align}\label{Gen}
\sL_t:=\tr((\sigma\sigma^*)(t,\bbx)\cdot\nabla^2_{x_1})+x_1\cdot\nabla_{x_2}+b(t,\bbx)\cdot \nabla_{x_1},
\end{align}
where $\tr$ denotes the trace of a matrix, and $\nabla_{x_i}$ denotes the gradient with respect to variable $x_i$.

The SDE \eqref{1}, also called stochastic Hamiltonian/kinetic system, 
represents essential frameworks for studying the dynamic behavior of physical, chemical, and biological systems under the influence of both deterministic laws and stochastic forces.
These systems integrate the principles of Hamiltonian mechanics with the inherent randomness captured by stochastic processes, providing a comprehensive understanding of complex phenomena across various disciplines. We refer the reader to \cite{So94,GSS78,K07} for a general background introduction, to \cite{GP14,HN04,Vi09} for the ergodicity of the system, and to \cite{MSH02,Ta02} for Euler approximations of the invariant measures. 

In recent years, there has been a great interest in the study of the well posedness for SDE \eqref{1} with irregular drifts and non-degenerate diffusion coefficients. 
For the strong well posedness of degenerate SDE \eqref{1}, Chaudru de Raynal in \cite{Ch17} obtained the first result under the assumptions that $b$ is $\alpha$-H\"older continuous in $x_1$ and $\beta$-H\"older continuous in $x_2$ with $\alpha\in(0,1)$ and $\beta\in(\frac{2}{3},1)$. This result was then extended by Wang and Zhang \cite{WZ16} under similar H\"older-Dini conditions. 
When $b$ is in some anisotropic Sobolev space, Zhang in \cite{Zh18} established the strong well-posedness of \eqref{1} by solving the associated Kolmogorov equation and 
Zvonkin's transformation (see also \cite{FFPV17} for similar result with constant diffusion coefficients but stronger drift condition than \cite{Zh18}).  
For the weak well-posedness of SDE \eqref{1}, to the best of our knowledge, the weakest result is obtained by Chaudru de Raynal and Menozzi in \cite[Theorem 1]{CM22}  under the assumptions that
$\sigma$ is bounded measurable and uniformly elliptic and H\"older continuous in $x$ uniformly in $t$, and
$$
b\in L^q(\mR_+; L^p(\mR^{2d}))\ \mbox{ with } \tfrac2q+\tfrac{4d}p<1.
$$
In fact, they considered more general models of SDEs  that can be seen as perturbed ODEs for which a noise acting on the first component will transmit to the whole chain
of ODEs through a weak H\"ormander like condition. It is noted that they also provide counterexamples in \cite[Theorem 2]{CM22} to illustrate the almost sharpness of condition $\tfrac2q+\tfrac{4d}p<1$.
We mention that kinetic SDEs driven by $\alpha$-stable processes were studied in \cite{HWZ21}.

Now we consider a second order singular interacting particle system in $\mR^d$:
\begin{align}\label{JH1}
\ddot X^{N,i}_t=\sum_{j\not=i}^N\frac{\gamma_j(X^{N,i}-X^{N,j})}{|X^{N,i}-X^{N,j}|^\alpha}-\gamma \dot X^{N,i}_t+\sqrt{2\nu}\dot W^{N,i}_t,\ \ i=1,\cdots,N,
\end{align}
where $\gamma_j\in\mR$, $\alpha\in[1,d+1)$, $\gamma\geq 0$ represents the damping force, and $\nu>0$ is the intensity of the noise, $(W^{N,i})_{i=1}^N$ is a sequence of i.i.d. Brownian motions.
Note that if we let $\dot X^{N,i}_t=:V^{N,i}_t$, then the above second order system can be written as the following first order system:
$$
\dot X^{N,i}_t=V^{N,i}_t,\ \ \dot V^{N,i}_t=\sum_{j\not=i}^N\frac{\gamma_j(X^{N,i}-X^{N,j})}{|X^{N,i}-X^{N,j}|^\alpha}-\gamma V^{N,i}_t+\sqrt{2\nu}\dot W^{N,i}_t.
$$
In particular, $Z^N_t:=(V^{N,i}_t, X^{N,i}_t)_{i=1,\cdots, N}$ satisfies a degenerate SDE in phase space $\mR^{2dN}$. 
It is easy to see that $x\mapsto x/|x|^\alpha \in L^p_{loc}(\mR^d)$ for any $p<d/(\alpha-1)$. Thus, to apply the result in \cite{CM22} to SDE \eqref{JH1}, 
one needs to assume $\alpha\in[1,1+1/(4N))$ due to the restriction $\tfrac2q+\tfrac{4Nd}p<1$. A natural question is that for a fixed $\alpha>1$, is it possible to establish the well-posedness of SDE \eqref{JH1}
for any $N\in\mN$? We shall give an affirmative answer that it is true for any $\alpha\in(1,\frac43)$ by considering the drift $b$ being in some Kato's class. In the non-degenerate case, this problem has been tackled in \cite{LX22, HRZ22}. Therein, the strong well-posedness are considered.

Another aim of this paper is to show the existence and two-sided estimates of the density or heat kernel of SDE \eqref{1} with singular drifts.
Let us first consider a special case $b=0$ and $\sigma=\mI_{d\times d}$:
\begin{align}
\dif X_t^1 =\sqrt{2} \dif W_t,\ \ \dif X_t^2 = X_t^1\dif t,\ \ t > 0.
\end{align}
In \cite{Kol34}, Kolmogorov first wrote down the density of $\bbX_t$ starting from $\bbX_0=\bbx\in\mR^{2d}$:
\begin{align}
p(0,\bbx;t,\bby)=\left(\frac{\sqrt{3}}{2\pi t^2}\right)^d\exp\left(-\frac{1}{4}\left|K_t^{-\frac{1}{2}}(\bby-\theta_t(\bbx))\right|^2\right),
\end{align}
where
\begin{align}
\theta_t(\bbx)=(x_1,x_2+x_1t),\ \ 
K_t=\left(
\begin{matrix}
2t\mI_{d\times d},&t^2\mI_{d\times d}\\
t^2\mI_{d\times d},&\frac{2t^3}{3}\mI_{d\times d}
\end{matrix}
\right).
\end{align}
From the expression, it is easy to derive the following sharp two-sided estimates:
\begin{align}\label{12}
\Big(\frac{\sqrt3}{2\pi t^2}\Big)^d\e^{-\frac{4+\sqrt{13}}{4}|\mT_t(\bby-\theta_t(\bbx))|^2}\leq p(0,\bbx;t,\bby)
\leq\Big(\frac{\sqrt3}{2\pi t^2}\Big)^d\e^{-\frac{4-\sqrt{13}}{4}|\mT_t(\bby-\theta_t(\bbx))|^2},
\end{align}
where for $t>0$, $\mT_t$ is the scale matrix given by
\begin{align}\label{215}
\mT_t=
\left(
\begin{matrix}
t^{-\frac{1}{2}}\mI_{d\times d} & 0_{d\times d}\\
0_{d\times d} & t^{-\frac{3}{2}}\mI_{d\times d}
\end{matrix}
\right).
\end{align}
Indeed, noting that
$$
K^{-1}_t=\begin{pmatrix}
2t^{-1} \mI_{d\times d}& -3t^{-2}\mI_{d\times d}\\ -3t^{-2}\mI_{d\times d} & 6t^{-3}\mI_{d\times d},
\end{pmatrix}
$$
the two-sided estimates \eqref{12} follow by
$$
-(2+\sqrt{13}) t|x_1|^2+(2-\sqrt{13})t^{-1}|x_2|^2\leq 6\<x_1, x_2\>\leq
-(2-\sqrt{13}) t|x_1|^2+(2+\sqrt{13})t^{-1}|x_2|^2.
$$

Since the seminal work of \cite{Kol34}, there are many works concentrated on the heat kernel estimates of SDE \eqref{1}. In \cite{KMM10}, the global upper and lower diagonal bounds for the heat kernel were obtained when the drift is uniformly Lipschitz continuous and bounded. The result was extended to more general models of SDEs  that can be seen a full chain of perturbed ODEs. We refer to the works \cite{DM10,Me11} for related discussions. The authors in \cite{CMPZ23} proved the two-sided heat kernel estimates of SDE \eqref{1} with unbounded and H\"older drifts. However, the drift term of the SDE discussed in \cite{KMM10,DM10,Me11,CMPZ23} does not include a singular drift.

Before stating our main result, we first introduce the following notion of weak solutions.
\bd\label{weak}
Let $\mathfrak{F}:=(\Omega, \sF, (\sF_t)_{t\ge0}, \bP)$ be a stochastic basis.
We call  a pair of $\sF_t$-adapted processes $(\bbX,W)$ on $\mathfrak{F}$ a weak solution of SDE \eqref{1} starting from point $\bbx\in\mR^{2d}$  if 
\begin{enumerate}[(i)]
\item $\bP(\bbX_0=\bbx)=1$;
\item $W$ is a $d$-dimensional $\sF_t$-Brownian motion;
\item For Lebesgue almost all $t>0$, the law of $\bbX_t$ is absolutely continuous with respect to the Lebesgue measure;
\item  For any $t>0$, it holds that $\int^t_0|b(s,\bbX_s)|\dif s<\infty$ a.s. and
\begin{align*}
\bbX_t=\bbx+\int^t_0 B(s,\bbX_s)\dif s+\int^t_0\Sigma(s,\bbX_s)\dif W_s,\ \ a.s.
\end{align*}
\end{enumerate}
\ed
\br
In the above definition, we require (iii) that is used to show the apriori Krylov estimate.
Without (iii), one can not use the mollifying technique to take limits (see Lemma \ref{Le41} below). It is an open question to show (iii) priorly under (i) (ii) and (iv). It is noted that in the non-degenerate case, we can drop it (see \cite{XXZZ}).
\er

Now we introduce the following main assumptions used below:
\begin{enumerate}[{\bf (H$^\sigma_b$)}]
\item There is a $\kappa_0\geq 1$ such that for all $t\geq 0$, $\bbx=(x_1,x_2)\in\mR^{2d}$ and $\xi\in\mR^d$,
\begin{align}
\kappa_0^{-1}|\xi|^2\le |\sigma(t,\bbx)\xi|^2\le \kappa_0|\xi|^2,
\end{align}
and for some $\gamma_0\in(0,1)$,
\begin{align}\label{31}
|\sigma(t,\bbx)-\sigma(t,\bby)|\le \kappa_0\big(|(\bbx-\bby)_1|^{\gamma_0}+|(\bbx-\bby)_2|^{\gamma_0/3}\big).
\end{align}
Moreover, $b=b_0+b_1$, where $b_0$ satisfies that for all $j=1,2,\cdots$,
\begin{align}\label{23}
|b_0(t,0)|\leq \kappa_0,\ \ \| \nabla^j b_0(t,\cdot)\|_\infty\leq \kappa_j, 
\end{align}
and $|b_1|$ belongs to Kato's class $\mK^{\rm o}_1$ (see Definition \ref{Def0} below).
\end{enumerate}
\br
If $|b_1|$ belongs to Kato's class $\mK^{\rm o}_1$ and $b_0$ only satisfies that for some $\beta\in(0,1]$,
$$
|b_0(t,0)|\leq\kappa_0, \ \ |b_0(t,\bbx)-b_0(t,\bby)|\leq\kappa_0|\bbx-\bby|^\beta,
$$
then $b$ still satisfies the assumption in {\bf (H$^\sigma_b$)}. In fact, we can let $\tilde b_0(t,\bbx):=b_0(t,\cdot)*\phi_n(\bbx)$ and $\tilde b_1:=b_1+b_0-\tilde b_0$, where $\{\phi_n\}_n$ is a family of modifiers. Then $\tilde b_0$ satisfies \eqref{23}, and $b_0-\tilde b_0$ is bounded measurable, and thus $\tilde b_1\in\mK^{\rm o}_1$.
\er

Let  $B_0(t,\bbx):=(b_0(t,\bbx),x_1)^*$ and $\theta_{t,s}(\bbx)$ be the unique solution of ODE
\begin{align}\label{100}
\dot{ \t}_{t,s}(\bbx)=B_0(t,\t_{t,s}(\bbx)),\ \ t, s\geq 0,\ \ { \t}_{s,s}(\bbx)=\bbx.
\end{align}
Note that $\{\theta_{t,s}(\bbx), t\geq s\geq 0\}$ stands for the forward flow and 
$\{\theta_{t,s}(\bbx), 0\leq t\leq s\}$ stands for the backward flow.
For $T\in(0,\infty)$ we introduce the following notation:
$$
\mD_T:=\{(s,\bbx;t,\bby): 0\leq s<t\leq T, \bbx,\bby\in\mR^{2d}\}.
$$

The aim of this paper is to prove the following theorem.
\bt\label{Main1}
Under {\bf (H$^\sigma_b$)}, for each $s>0$ and $\bbx\in\mR^{2d}$, there is a unique weak solution $(\bX_{t,s}(\bbx))_{t\geq s}$ to SDE \eqref{1} starting from  $\bbx$ at time $s$ in the sense of Definition \eqref{weak}.
Moreover, if $|b_1|^\g\in \mK_1$ for some $\gamma>1$, then $\bX_{t,s}(\bbx)$ admits a density $p(s,\bbx;t,\bby)$ enjoying the following two-sided estimates: for any $T>0$ and all $(s,\bbx;t,\bby)\in\mD_T$,
\begin{align}\label{AX1}
C_0(t-s)^{-2d}\e^{-\lambda_0|\mT_{t-s}(\t_{t,s}(\bbx)-\bby)|^2}\le p(s,\bbx;t,\bby) \le C_1(t-s)^{-2d}\e^{-\lambda_1|\mT_{t-s}(\t_{t,s}(\bbx)-\bby)|^2},
\end{align}
where $\lambda_0,\lambda_1$ depend on $T, d, \gamma_0, \kappa_0, \kappa_1$, and $C_0,C_1$ depend on $T, d, \gamma_0, \kappa_0, \kappa_1$ and function $|b_1|^\g$.
\et

\br
By Remark \ref{Re24} below and the above theorem, one sees that for any $\alpha\in(1,\frac43)$, SDE \eqref{JH1} admits a unique weak solution and a density for any starting point $\bbx$.
\er

This paper is organized as follows: In Section 2, we give the definition of Kato's class in kinetic case and prove a crucial estimate  in  Lemma \ref{Le26} related to the function in Kato's class.
Roughly to say, let $p_{s,t}(x,\cdot)$
be a family of time-inhomoegenous transition kernels. For a function $b:\mR_+\times\mR^d\to [0,\infty)$, we show an estimate  in  Lemma \ref{Le26} as
$$
\int^t_s\!\!\!\int_{\mR^d} p_{s,r}(x,\dif z) b(r,z) p_{r,t}(z,A)\dif r\leq K_b\cdot p_{s,t}(x,A),
$$
where $K_b$ is some functional of $b$ that is used to define the Kato class. Such an estimate plays a crucial role in establishing the two-sided estimates of heat kernel.

In Section 3, we  solve the kinetic PDEs with drifts being unbounded and in Kato's class. The main result is Theorem \ref{37}. To figure out the main idea of the proof, we consider the following kinetic SDE:
$$
\p_t u+\Delta_{x_1} u+x_1\cdot\nabla_{x_2}u+b\cdot\nabla_{x_1}u+f=0,\ \ u(T)=0,
$$
where $b=b_0+b_1$ satisfies the assumptions in {\bf (H$^\sigma_b$)}.
Since $b_0$ is an unbounded vector field, one can not directly solve it by considering the following integral equation
$$
u(s,\bbx)=\int^T_s P_{r-s}(b\cdot\nabla_{x_1}u+f)(s,\bbx)\dif r,
$$
where $P_t$ is the kinetic semigroup associated with $\Delta_{x_1} +x_1\cdot\nabla_{x_2}$.
Instead of this, we consider the following perturbation integral equation
\begin{align}\label{int}
u(s,\bbx)=\int^T_s P^{(0)}_{r,s}(b_1\cdot \nabla_{x_1} u+f)(r,\bbx)\dif r,
\end{align}
where $P^{(0)}_{t,s}\varphi(\bbx):=\mE\varphi(X_{t,s}(\bbx))$ and $X_{t,s}(\bbx)$ is the unique solution of SDE \eqref{1} with $\sigma=\mI$ and $b=b_0$, 
and starting from $\bbx$ at time $s$. To solve \eqref{int}, we shall use the heat kernel estimates obtained in \cite{CMPZ23}.

In Section 4, we establish the weak well-posedness of SDE \eqref{1} under {\bf (H$^\sigma_b$)} based on the solvability of the kinetic PDE. In Section 5, we prove the two-sided estimates of the density following the same line as in \cite{CMPZ23}.
However, to treat the singular part, we need a Krylov-type estimate for the controlled SDE (see Lemma \ref{Kry} below).

We conclude this introduction by introducing the following convention: 
throughout this paper, we use $C$ with or without subscripts to denote an unimportant constant, whose value
may change in different occasions. We also use $:=$ as a way of definition. By $A\lesssim_C B$ or simply $A\lesssim B$, we mean that $A\leq C B$, for some unimportant constant $C\geq 1$.

\section{Preliminaries}

The Kato class was first introduced by Kato in \cite{Ka72} to prove the essential self-adjointness of Schrödinger operator. It plays important roles in the analysis of partial differential equations and mathematical physics. A measurable function $f$ on $\mR^d$ is said to be of Kato class $\mK_d$ if
\begin{align*}
\lim_{r\to 0}\sup_{x\in\mR^d}\int_{|x-y|\le r}|f(y)| K(x-y) \dif y=0,
\end{align*}
where 
$$
K(x)=|x|^{2-d}\cdot \1_{\{d\ge 3\}}+ \log|x|^{-1}\cdot \1_{\{d=2\}}+|x|\cdot \1_{\{d=1\}}.
$$
In \cite{AS82}, Aizenman and Simon provided an equivalent characterization for $\mK_d$ by
\begin{align*}
\lim_{\delta\to 0}\sup_{x\in\mR^d}\int_0^\delta \int_{\mR^d} (2\pi s)^{-\frac d 2}\e^{\frac{|x-y|^2}{2s}}|f(y)| \dif y \dif s=0.
\end{align*}
We refer to \cite[Theorem 3.6]{CZ95} for more information of this result. For general time-dependent case, Q.S.Zhang in \cite{ZQ96} introduce the following definition of Kato’s class:
\begin{align}\label{K0}
\lim_{\delta\to 0}\sup_{t>0}\sup_{x\in\mR^d}\int_0^\delta \int_{\mR^d} s^{-\frac {d+1} 2}\e^{\frac{\lambda|y|^2}{s}}|f(t\pm s,x\pm y)| \dif y \dif s=0, \qquad \forall \lambda>0,
\end{align}
and showed two-sided Gaussian estimates for second order operators with drift satisfying \eqref{K0} in \cite{ZQ97}.

In this section we introduce the Kato class of space-time functions on $\mR\times\mR^{2d}$, 
which extends  the classical definition of Kato’s class in \cite{ZQ96} to kinetic case. Then we present some basic properties of functions in Kato's class
used in this paper. Moreover, we show a crucial estimate in Lemma \ref{Le26} that shall be used to prove the two-sided estimates of heat kernel.

To introduce the Kato class, we need the following Gaussian kernel function as used in \eqref{12}: for $\beta\geq 0$ and $\lambda>0$,
\begin{align}\label{24}
\eta^{(\beta)}_\lambda(t,\bbx):=t^{-\frac\beta2-2d}\e^{-\lambda |\mT_t\bbx|^2},\ \ (t,\bbx)\in\mR_+\times\mR^{2d},
\end{align}
where $\mT_t$ is defined by \eqref{215}.
Clearly,
\begin{align}\label{214}
\eta^{(\beta)}_\lambda(1,\cdot)\in\cap_{p\in[1,\infty]}L^p(\mR^{2d}).
\end{align}
For $\lambda,\delta>0$ and a space-time function $f(t,\bbx):  \mR^{1+2d}\to\mR$, we define
\begin{align}\label{Def1}
K^{(\beta)}_\lambda(f; \delta)
:=\sup_{t\geq 0}\sup_{\bbx\in\mR^{2d}}\bigg\{\int_0^\delta\!\!\int_{\mR^{2d}}\eta^{(\beta)}_\lambda(r,\bby)|f(t\pm r,\bbx\pm\bby)|\dif \bby \dif r\bigg\},
\end{align}
where
$$
|f(t\pm r,\bbx\pm\bby)|:=\sum_{i,j=0,1}|f(t+(-1)^ir,\bbx+(-1)^j\bby)|.
$$
It is easy to see that $\delta\mapsto K^{(\beta)}_\lambda(f; \delta)$ is monotone increasing and
$\lambda\mapsto K^{(\beta)}_\lambda(f; \delta)$ is monotone decreasing, and for $0\leq \beta'\leq\beta$,
\begin{align}\label{AA1}
K^{(\beta')}_\lambda(f; \delta)\leq \delta^{\frac{\beta-\beta'}{2}}K^{(\beta)}_\lambda(f; \delta),\ \ \forall\delta>0. 
\end{align}

We introduce an anisotropic distance in $\mR^{2d}$ by
$$
|\bbx|_{\rm d}:=\max\{|x_1|,|x_2|^{\frac 13}\}, \qquad \bbx=(x_1,x_2)\in\R^{2d}.
$$
The following lemma tells us that for $\beta,\lambda>0$, 
$$
\mbox{
$K^{(\beta)}_\lambda(f;1)<\infty\Rightarrow f\in L^1_{loc}(\mR^{1+2d})$.}
$$
\bl\label{Le23}
For any $\beta>0$, there is a constant $C=C(d,\beta)>0$ such that for all  $\lambda,\delta\in(0,1)$,
$$
\sup_{(t,\bbx)\in\mR^{1+2d}}\int^{\delta}_0\!\!\!\int_{|\bby|_{\rm d}\leq\sqrt\delta}|f(t\pm r,\bbx\pm\bby)|\dif \bby \dif r\leq  C\e^{10\lambda}\delta^{\frac\beta2+2d}K^{(\beta)}_\lambda(f; \delta).
$$
\el
\begin{proof}
By definition we have
\begin{align*}
&\int_0^{\delta}\!\!\!\int_{\mR^{2d}}\eta^{(\beta)}_\lambda(r,\bby)|f(t\pm r,\bbx\pm\bby)|\dif \bby \dif r\\
&\quad\geq\int^{\delta}_{\delta/2}\!\int_{|\bby|_{\rm d}\leq\sqrt\delta} r^{-\frac\beta2-2d}\e^{-\lambda |\mT_r\bby|^2}|f(t\pm r,\bbx\pm\bby)|\dif \bby \dif r\\
&\quad\geq \delta^{-\frac\beta2-2d}\e^{-10\lambda}\int^{\delta}_{\delta/2}\!\int_{|\bby|_{\rm d}\leq\sqrt\delta}|f(t\pm r,\bbx\pm\bby)|\dif \bby \dif r.
\end{align*}
Hence,
$$
\sup_{(t,\bbx)\in\mR^{1+2d}}\int^{\delta}_{\delta/2}\!\int_{|\bby|_{\rm d}\leq\sqrt\delta}|f(t\pm r,\bbx\pm\bby)|\dif \bby \dif r\leq \e^{10\lambda}\delta^{\frac\beta2+2d}K^{(\beta)}_\lambda(f; \delta).
$$
On the other hand, in view that  there exist $M_n:=d^{2d}2^{2nd}$ points  $\{\bbx_i\}\subset\mR^{2d}$ so that
$$
\big\{\bby: |\bby|_{\rm d}\leq\sqrt\delta\big\}\subset \cup_{i=1}^{M_n}\big\{\bby: |\bby-\bbx_i|_{\rm d}\leq\sqrt{2^{-n}\delta}\big\},
$$
by the change of variable, we have
\begin{align*}
&\int^{\delta}_0\!\!\!\int_{|\bby|_{\rm d}\leq\sqrt\delta}|f(t\pm r,\bbx\pm\bby)|\dif \bby \dif r
=\sum_{n=0}^\infty\int^{2^{-n}\delta}_{2^{-(n+1)}\delta}\int_{|\bby|_{\rm d}\leq\sqrt\delta}|f(t\pm r,\bbx\pm\bby)|\dif \bby \dif r\\
&\qquad\leq\sum_{n=0}^\infty\sum_{i=1}^{M_n}\int^{2^{-n}\delta}_{2^{-(n+1)}\delta}\int_{|\bby-\bbx_i|_{\rm d}\leq\sqrt{2^{-n}\delta}}|f(t\pm r,\bbx\pm\bby)|\dif \bby \dif r\\
&\qquad\leq\sum_{n=0}^\infty M_n\sup_i\int^{2^{-n}\delta}_{2^{-(n+1)}\delta}\int_{|\bby|_{\rm d}\leq\sqrt{2^{-n}\delta}}|f(t\pm r,\bbx\pm(\bbx_i+\bby))|\dif \bby \dif r\\
&\qquad\leq\sum_{n=0}^\infty d^{2d}2^{2nd}\e^{10\lambda}(2^{-n}\delta)^{\frac\beta2+2d}K^{(\beta)}_\lambda(f; 2^{-n}\delta)
=\frac{d^{2d}\delta^{\frac\beta2+2d}\e^{10\lambda}}{1- 2^{-\beta/2}}K^{(\beta)}_\lambda(f; 2^{-n}\delta).
\end{align*}
The proof is complete.
\end{proof}

Now we introduce the definition of Kato's class in the kinetic case used below.
\bd[Kato's class]  \label{Def0}
For any $\beta\geq 0$, we define
\begin{align*}
\mK_{\beta}&:=\Big\{f\in L^1_{loc}(\mR^{1+2d}):  K^{(\beta)}_\lambda(f;1)<\infty,\mbox{ $\forall \lambda>0$} \Big\},\\
\mK^{\rm o}_{\beta}&:=\bigg\{f\in \mK_{\beta}: \lim_{\delta\downarrow 0}K^{(\beta)}_\lambda(f;\delta)=0,\mbox{ $\forall \lambda>0$} \bigg\},
\end{align*}
and
\begin{align*}
\mK^{\rm p}_{\beta}:=\bigg\{f\in\mK^{\rm o}_\beta:  \mbox{  $\forall p,\lambda>0$, $\exists N\in\mN$  s.t. } \int_0^1\frac{(K^{(\beta)}_\lambda(f;r))^{N}}{r^p}\dif r<\infty \bigg\}.
\end{align*}
\ed
\br
 If $f(t,\bbx)=f(t,x_1)$ only depends on the first variable $x_1$, then the definition of $\mK^{\rm o}_{1}$ coincides with the classical definition \eqref{K0}.
By \eqref{AA1}, one sees that for $0\leq \beta'<\beta$,
$$
\mK_\beta\subset\mK^{\rm p}_{\beta'}\subset \mK^{\rm o}_{\beta'}.
$$
For a function $f: [0,T]\times\mR^{2d}\to\mR$, we shall extend it to $\mR\times\mR^{2d}$ by putting $f(t,\cdot)=0$ for $t\notin[0,T]$.
\er

\br\label{Re24}
Let $f(x)=|x_2|^{1-\alpha}$. It is easy to see that for $\alpha\in(1,\frac43)$, $f\in\mK^{\rm o}_1$. Indeed, by definition and the change of variable, we have
\begin{align*}
K^{(1)}_\lambda(f; \delta)
&=\sup_{\bbx\in\mR^{2d}} \int_0^\delta\!\!\int_{\mR^{2d}}r^{-\frac12-2d}\e^{-\lambda |\mT_r\bby|^2}|x_2-y_2|^{1-\alpha}\dif \bby \dif r \\
&=\sup_{x_2\in\mR^{d}}\int_0^\delta\!\!\int_{\mR^{d}}r^{-\frac{1+3d}2}\e^{-\lambda |y_2|^2/r^3}|x_2-y_2|^{1-\alpha}\dif y_2 \dif r\\
&=\sup_{x_2\in\mR^{d}}\int_0^\delta\!\!\int_{\mR^{d}}r^{1-\frac{3\alpha}2}\e^{-\lambda |x_2-y_2|^2}|y_2|^{1-\alpha}\dif y_2 \dif r\\
&\leq \frac{2\delta^{2-\frac{3\alpha}2}}{4-3\alpha}\left(\int_{|y|\leq 1}|y|^{1-\alpha}\dif y+\int_{\mR^d}\e^{-\lambda |y|^2}\dif y\right).
\end{align*}
\er

For $n\in\mN$ and $f\in\mK_\beta$, let 
$$
f_n(t,\bbx):=(f*\phi_n)(t,\bbx)=\int_\mR\int_{\mR^{2d}}f(t-s,\bbx-\bby)\phi_n(s,\bby)\dif\bby\dif s,
$$
where for a probability density function $\phi\in C^\infty_c(\mR^{1+2d})$,
\begin{align}\label{MO1}
\phi_n(s,\bbx):=n^{1+2d}\phi(ns, n\bbx).
\end{align}
By the definition of convolution, it is easy to see that
\begin{align}\label{CC1}
K^{(\beta)}_\lambda(f_n; \delta)=K^{(\beta)}_\lambda(|f_n|; \delta)\leq K^{(\beta)}_\lambda(|f|; \delta)=K^{(\beta)}_\lambda(f; \delta).
\end{align}
Moreover, if $\beta>0$, then by Lemma \ref{Le23}, one sees that
$$
f_n\in C^\infty_b(\mR^{1+2d}).
$$

The following lemma is useful.
\bl\label{Le24}
Let $\beta\geq 0$. For any $f\in \mK^{\rm o}_\beta$ and $T,R>0$, it holds that
$$
\lim_{n\to\infty}K^{(\beta)}_\lambda(\1_{D_R}|f_n-f|;T)=0,
$$
where $f_n=f*\phi_n$ and $D_R:=\{(t,\bbx)\in\mR^{1+2d}: |t|\leq R,\ |\bbx|_{\rm d}\leq R\}$.
\el
\begin{proof}
For simplicity of notations, we write $g^R_n:=\1_{D_R}|f_n-f|$.
For any $\delta\in(0,T)$, we have
\begin{align*}
K^{(\beta)}_\lambda(g^R_n; T)
&\leq K^{(\beta)}_\lambda(g^R_n; \delta)
+\sup_{t\geq 0}\sup_{\bbx\in\mR^{2d}}\bigg\{\int^T_\delta\!\!\int_{\mR^{2d}}\eta^{(\beta)}_\lambda(r,\bby)g^R_n(t\pm r,\bbx\pm\bby)\dif \bby \dif r\bigg\}\\
&\leq K^{(\beta)}_\lambda(g^R_n; \delta)
+\delta^{-\frac\beta2-2d}\sup_{t\geq 0}\sup_{\bbx\in\mR^{2d}}\bigg\{\int^T_\delta\!\!\int_{\mR^{2d}}g^R_n(t\pm r,\bbx\pm\bby)\dif \bby \dif r\bigg\}\\
&\leq 2K^{(\beta)}_\lambda(|f|; \delta)
+\delta^{-\frac\beta2-2d}\int_{\mR^{1+2d}}g^R_n(r,\bby)\dif \bby \dif r\\
&=2K^{(\beta)}_\lambda(|f|; \delta)
+\delta^{-\frac\beta2-2d}\int_{D_R}|f_n-f|(r,\bby)\dif \bby \dif r,
\end{align*}
which converges to zero by first letting $n\to\infty$ and then $\delta\to 0$.
\end{proof}

The following lemma asserts that the mixed $L^q_tL^\bbp_\bbx$-space  is contained in $\mK^{\rm p}_\beta$.
\bl
Let $q\in[1,\infty]$ and $\bbp=(p_1,\cdots,p_{2d})\in[1,\infty]^{2d}$ satisfy 
$$
\kappa:=2-(\tfrac{1}{p_1}+\cdots\tfrac1{p_d}+\tfrac3{p_{d+1}}+\cdots+\tfrac{3}{p_{2d}}+\tfrac{2}{q})>0.
$$ 
For any $\beta\in[0,\kappa)$,
there is a constant $C=C(d,\bbp,q,\beta)>0$ such that for all $T\in[0,1]$,
$$
K^{(\beta)}_\lambda(f;T)\lesssim_CT^{(\kappa-\beta)/2}\|f\|_{\mL^q(\mL^{\bbp})},
$$
where 
$$
\|f\|_{\mL^q(\mL^{\bbp})}=\left(\int_\mR\left(\int_{\mR}\cdots\left(\int_{\mR}|f(s,x^{(1)},\cdots,x^{(2d)})|^{p_1}\dif x^{(1)}\right)^\frac{p_2}{p_1}\cdots\dif x^{(2d)}\right)^\frac{q}{p_{2d}}\dif s\right)^\frac{1}{q}.
$$
In particular, for any $\beta\in[0,\kappa)$,
$$
\mL^q(\mL^{\bbp})\subset\mK^{\rm p}_\beta.
$$
\el
\begin{proof}
Let $\bar q:=\frac{q}{q-1}$ and $\bar\bbp=(\bar p_1,\cdots, \bar p_{2d})$ with $\bar p_i:=\frac{p_i}{p_i-1}$.
By H\"older's inequality, we have
$$
\int_0^T\!\!\!\int_{\mR^{2d}}\eta^{(\beta)}_\lambda(s,\bby)|f(t+s, \bbx+\bby)|\dif \bby\dif s\le \|\eta^{(\beta)}_\lambda\|_{\mL^{\bar q}_T(\mL^{\bar\bbp})}
\|f\|_{\mL^{q}(\mL^{\bbp})},
$$
where
$$
\|\eta^{(\beta)}_\lambda\|_{\mL_T^{\bar q}(\mL^{\bar\bbp})}=\left(\int_0^T\left(\int_{\mR}\cdots\left(\int_{\mR}|\eta^{(\beta)}_\lambda(s,x^{(1)},\cdots,x^{(2d)})|^{\bar p_1}\dif x^{(1)}\right)^\frac{\bar p_2}{\bar p_1}\cdots\dif x^{(2d)}\right)^\frac{\bar q}{\bar p_{2d}}\dif s\right)^\frac{1}{\bar q}.
$$
By the change of variables, it is easy to see that for each $s\in(0,1)$,
$$
\|\eta^{(\beta)}_\lambda(s)\|_{\bar\bbp}=s^{\gamma}\|\eta^{(\beta)}_\lambda(1)\|_{\bar\bbp},
$$
where $\eta^{(\beta)}_\lambda(s)=\eta^{(\beta)}_\lambda(s,\cdot)$ and
$$
\gamma:=\tfrac{1}{2\bar p_1}+\cdots+\tfrac{1}{2\bar p_d}+\tfrac{3}{2\bar p_{d+1}}+\cdots+\tfrac{3}{2\bar p_{2d}}-2d-\tfrac\beta2.
$$
Hence,
$$
\|\eta^{(\beta)}_\lambda\|_{\mL^{\bar q}_T(\mL^{\bar\bbp})}=\left(\int_0^Ts^{\gamma\bar q}
\dif s\right)^\frac{1}{\bar q}\|\eta^{(\beta)}_\lambda(1)\|_{\bar\bbp}
\stackrel{\eqref{214}}{\lesssim} T^{\gamma+\frac{1}{\bar q}}=T^{\frac{\kappa-\beta}2}.
$$
The proof is complete.
\end{proof}
We also need the following lemma.
\bl\label{KC}
For $\beta\in[0,2)$ and $\g>1$, if $|f|^\g\in\mK_\beta$, then $f\in\mK^{\rm p}_\beta$.
\el
\begin{proof}
Let $\beta\in[0,2)$, $\gamma>1$ and $\lambda>0$.
By H\"older's inequality, we have for any $r\in(0,1)$,
\begin{align*}
K^{(\beta)}_\lambda(|f|; r)
&\le K^{(\beta)}_\lambda(|f|^\g; 1)\left(\int_0^r\!\!\!\int_{\mR^{2d}}\eta^{(\beta)}_\lambda(s,\bby)\dif\bby\dif s\right)^\frac{\g-1}{\g}\\
&=K^{(\beta)}_\lambda(|f|^\g; 1)\left(\int_0^r s^{-\frac{\beta}{2}}\dif s\int_{\mR^{2d}}\eta^{(\beta)}_\lambda(1,\bby)\dif\bby\right)^\frac{\g-1}{\g}\\
& = K^{(\beta)}_\lambda(|f|^\g; 1)\left(\int_{\mR^{2d}}\eta^{(\beta)}_\lambda(1,\bby)\dif\bby\right)^\frac{\g-1}{\g}\Big(\tfrac{2}{2-\beta}r^{\frac{2-\beta}{2}}\Big)^\frac{\g-1}{\g}.
\end{align*}
Thus, for any $p\geq 1$, one can choose $N>\frac{2\gamma(p-1)}{(2-\beta)(\g-1)}$ so that
$$
\int_0^1\frac{(K^{(\beta)}_\lambda(|f|; r))^{N}}{r^p}\dif r<\infty.
$$
The proof is complete.
\end{proof}

For $\beta\in\mR$ and $\l>0$, we introduce a basic function for later use:
\begin{align}\label{107}
g^{(\beta)}_{\l}(s,\bbx;t,\bby):=\eta^{(\beta)}_\lambda(t-s,\t_{t,s}(\bbx)-\bby),
\end{align}
where $\t_{t,s}(\bbx)$ is defined in \eqref{100}. Note that for $s\leq r\leq t$,
$$
\theta_{t,s}(\bbx)=\theta_{t,r}\circ \theta_{r,s}(\bbx).
$$

We recall the following estimates from  \cite[Lemma 2.2]{CMPZ23}:
\bl
Under \eqref{23}, there exists a constant $c_0=c_0(\k_1,d)\ge 1$ such that  for all $(s,\bbx;t,\bby)\in\mD_1$ and $r\in[s,t]$,
\begin{align}\label{16}
c_0^{-1}(|\mT_{t-s}(\bbx-\t_{r,t}(\bby))|-1)\le |\mT_{t-s}(\t_{t,r}(\bbx)-\bby)|\le c_0(|\mT_{t-s}(\bbx-\t_{r,t}(\bby))|+1).
\end{align}
\el

The following estimate plays a crucial role in the two-sided estimates of heat kernel.
\bl\label{Le26}
For any $\lambda>0,\beta\ge\alpha\ge 0$ and $0\leq f\in\mK_{\beta}$, there are constants $C_1=C_1(d,\alpha,\lambda)>0$ 
and $\kappa\in(0,1/2)$ depending only on $c_0$, where $c_0$ is from \eqref{16}, such that for all $(s,\bbx;t,\bby)\in\mD_1$,
\begin{align}\label{7}
\int_s^t\!\!\!\int_{\mR^{2d}}g^{(\alpha)}_{\lambda}(s,\bbx;r,\bbz) f(r,\bbz) g^{(\beta)}_{\lambda}(r,\bbz;t,\bby)\dif \bbz \dif r
\le C_1K^{(\beta)}_{\kappa\lambda}(f; t-s)g^{(\alpha)}_{\kappa\lambda}(s,\bbx;t,\bby).
\end{align}
\el
\begin{proof}
By  definitions \eqref{107} and \eqref{24}, since $f\geq 0$, we have
$$
\int_s^t\!\!\!\int_{\mR^{2d}}g^{(\alpha)}_{\lambda}(s,\bbx;r,\bbz)f(r,\bbz)g^{(\beta)}_{\lambda}(r,\bbz;t,\bby)\dif \bbz \dif r\leq I_1+I_2,
$$
where
\begin{align*}
I_1&=\int_s^{\frac{t+s}{2}}\!\!\!\!\int_{\mR^{2d}}(r-s)^{-\frac\alpha2-2d}(t-r)^{-\frac\beta2-2d}h_\lambda(s,x;r,z;t,y)f(r,\bbz)\dif \bbz \dif r,\\ 
I_2&=\int_{\frac{t+s}{2}}^t\!\!\int_{\mR^{2d}}(r-s)^{-\frac\alpha2-2d}(t-r)^{-\frac\beta2-2d}h_\lambda(s,x;r,z;t,y) f(r,\bbz)\dif \bbz \dif r,
\end{align*}
and
$$
h_\lambda(s,x;r,z;t,y):=\e^{-\l|\mT_{r-s}(\t_{r,s}(\bbx)-\bbz)|^2-\l|\mT_{t-r}(\t_{t,r}(\bbz)-\bby)|^2}.
$$
Observe that by  $(a+b)^2\le 2(a^2+b^2)$,
\begin{align}\label{E2}
h_\lambda(s,x;r,z;t,y)&=\e^{-\l(|\mT_{r-s}(\t_{r,s}(\bbx)-\bbz)|^2+|\mT_{t-r}(\t_{t,r}(\bbz)-\bby)|^2)}\no\\
&\lesssim \e^{-\l(|\mT_{r-s}(\t_{r,s}(\bbx)-\bbz)|^2+|\mT_{t-r}(\bbz-\t_{r,t}(\bby))|^2)/(2\c_0^2)}\nonumber\\
&\leq \e^{-\l(|\mT_{t-s}(\t_{r,s}(\bbx)-\bbz)|^2+|\mT_{t-s}(\bbz-\t_{r,t}(\bby))|^2)/(2\c_0^2)}\nonumber\\
&\leq \e^{-\l(|\mT_{t-s}(\t_{r,s}(\bbx)-\t_{r,t}(\bby))|^2)/(4\c_0^2)}\nonumber\\
&\stackrel{\eqref{16}}{\lesssim} \e^{-\l(|\mT_{t-s}(\t_{t,s}(\bbx)-\bby)|^2)/(8\c_0^4)}.
\end{align}
Let us first treat the term $I_1$. For $r\in(s,\tfrac{t+s}{2})$,
noting that  by $\beta\ge\alpha\ge 0$,
$$
(r-s)^{-\frac\alpha2}(t-r)^{-\frac\beta2}\leq 2^{\frac\alpha2}(r-s)^{-\frac\beta2}(t-s)^{-\frac\alpha2}
$$
and
$$
h_\lambda(s,x;r,z;t,y)\leq \e^{-\frac\l 2|\mT_{r-s}(\t_{r,s}(\bbx)-\bbz)|^2}h_{\lambda/2}(s,x;r,z;t,y),
$$
by definition \eqref{107}, for $\kappa:=1/(16\c_0^4)$, we have
\begin{align*}
I_1&\lesssim (t-s)^{-\frac\alpha2-2d}\int_s^{\frac{t+s}{2}}\!\!\!\int_{\mR^{2d}}g^{(\beta)}_{\l/2}(s,\bbx;r,\bbz)h_{\lambda/2}(s,x;r,z;t,y) f(r,z)\dif z\dif r\\
&\stackrel{\eqref{E2}}{\lesssim} (t-s)^{-\frac\alpha2-2d}\e^{-\kappa\l(|\mT_{t-s}(\t_{t,s}(\bbx)-\bby)|^2)}\int_s^t\!\!\!\int_{\mR^{2d}}g^{(\beta)}_{\l/2}(s,\bbx;r,\bbz)f(r,z)\dif z\dif r\\
&\le g^{(\alpha)}_{\kappa\l}(s,\bbx;t,\bby)\sup_{\bbx}\int_0^{t-s}\!\!\!\int_{\mR^{2d}}\eta^{(\beta)}_{\l/2}(r,\bbz)f(r+s,\bbz+\bbx)\dif\bbz\dif r\\
&\leq g^{(\alpha)}_{\kappa\l}(s,\bbx;t,\bby)K^{(\beta)}_{\l/2}(f; t-s).
\end{align*}
Similarly, for $I_2$, noting that
$$
h_\lambda(s,x;r,z;t,y)\leq \e^{-\frac\l 2|\mT_{t-r}(\t_{t,r}(\bbz)-\bby)|^2}h_{\lambda/2}(s,x;r,z;t,y),
$$
by definition, for $\kappa:=1/(16c_0^4)$, we have
\begin{align*}
I_2&\lesssim(t-s)^{-\frac\alpha2-2d}\int_{\frac{t+s}{2}}^t\!\!\int_{\mR^{2d}}g^{(\beta)}_{\l/2}(r,\bbz;t,\bby)h_{\lambda/2}(s,x;r,z;t,y)f(r,\bbz)\dif\bbz\dif r\\
&\stackrel{\eqref{E2}}{\lesssim} (t-s)^{-\frac\alpha2-2d}\e^{-\kappa\l(|\mT_{t-s}(\t_{t,s}(\bbx)-\bby)|^2)}\int_s^t\!\!\int_{\mR^{2d}}g^{(\beta)}_{\l/2}(r,\bbz;t,\bby)f(r,z)\dif z\dif r\\
&= g^{(\alpha)}_{\kappa\l}(s,\bbx;t,\bby)\int_s^t\!\!\int_{\mR^{2d}}\eta^{(\beta)}_{\l/2}(t-r,\t_{t,r}(\bbz)-\bby)f(r,z)\dif z\dif r\\
&\stackrel{\eqref{16}}{\lesssim} g^{(\alpha)}_{\kappa\l}(s,\bbx;t,\bby)\int_s^t\!\!\int_{\mR^{2d}}\eta^{(\beta)}_{\kappa\l}(t-r,\bbz-\t_{r,t}(\bby))f(r,z)\dif z\dif r\\
&\le g^{(\alpha)}_{\kappa\l}(s,\bbx;t,\bby)\sup_{\bby}\int_0^{t-s}\!\!\!\int_{\mR^{2d}}\eta^{(\beta)}_{\kappa\l}(r,\bbz)f(t-r,\bbz+\bby)\dif\bbz\dif r\\
&\leq g^{(\alpha)}_{\kappa\l}(s,\bbx;t,\bby)K^{(\beta)}_{\kappa\l}(f; t-s).
\end{align*}
The proof is complete.
\end{proof}

\section{Solvability of kinetic PDEs with drifts in Kato class}

In this section, we  assume the condition {\bf (H$^\sigma_b$)}. Firstly, for each $\bbx=(x_1,x_2)\in\mR^{2d}$, we consider the following SDE with regular drift coefficient $b_0$:
\begin{align}\label{SDE0}
\dif\bbX_{t,s}=B_0(t,\bbX_{t,s})\dif t+\Sigma(t,\bbX_{t,s})\dif W_t,\ \ \bbX_{s,s}=\bbx,\ \ t>s,
\end{align}
where $\Sigma$ is defined in \eqref{S} and
$$
B_0(t,\bbx):=(b_0(t,\bbx),x_1)^*.
$$
Under {\bf (H$^\sigma_b$)}, it is well-known that there exists a unique weak solution $\bbX_{t,s}(\bbx)$ to SDE \eqref{SDE0} (cf. \cite{CM22}).
For any $\varphi\in C_b(\mR^{2d})$ and $\bbx\in\mR^{2d}$, we define
\begin{align}\label{205}
P^{(0)}_{t,s}\varphi(\bbx):=\bE\varphi(\bbX_{t,s}(\bbx)).
\end{align}
Below we shall work on the time interval $[0,1]$ and use the following parameter set:
\begin{align}\label{Para}
\Theta:=(d,\gamma_0,\kappa_0,\kappa_1),
\end{align}
where $\gamma_0,\kappa_0,\kappa_1$ come from {\bf (H$^\sigma_b$)}.

We recall the following result proven in \cite[Theorem 1.1]{CMPZ23} about the density of $\bbX_{t,s}(\bbx)$.
\bt\label{Th31}
For each $0\leq s<t\leq 1$ and $\bbx\in\mR^{2d}$, $X_{t,s}(\bbx)$  admits a density $p_0(s,\bbx;t,\bby)$, which is continuous on $\mD_1$, so that, for any $\varphi\in C_b(\mR^{2d})$,
$$
P^{(0)}_{t,s} \varphi(\bbx)=\int_{\mR^{2d}}\varphi(\bby)p_0(s,\bbx;t,\bby)\dif\bby. 
$$
Moreover,  $p_0(s,\bbx;t,\bby)$ enjoys the following estimates:
\begin{enumerate}[(i)]
\item {(Two-sided estimates)} There are $\lambda\in(0,1)$ and $C_0\ge 1$ such that for all $(s,\bbx;t,\bby)\in\mD_1$,
\begin{align}\label{6}
C_0^{-1}g^{(0)}_{\lambda^{-1}}(s,\bbx;t,\bby)\le p_0(s,\bbx;t,\bby) \le C_0g^{(0)}_{\lambda}(s,\bbx;t,\bby).
\end{align}
\item{(Derivative estimate in $x_1$)} For $j=1,2$, there exist $\lambda,C_j>0$ such that  for all $(s,\bbx;t,\bby)\in\mD_1$,
\begin{align}\label{5}
|\nabla_{x_1}^jp_0(s,\bbx;t,\bby)|\le C_jg^{(j)}_{\lambda}(s,\bbx;t,\bby).
\end{align}
\item{(H\"older estimate in $\bbx$)} Let $\g_0,\g_1\in(0,1)$. For $j=0,1$, there exist constants $\lambda,C_j>0$ such that  for all $(s,\bbx;t,\bby)\in\mD_1$ and $\bbx'\in\mR^{2d}$,
\begin{align}\label{A1}
\begin{split}
&|\nabla_{x_1}^jp_0(s,\bbx;t,\bby)-\nabla_{x_1}^jp_0(s,\bbx';t,\bby)|\\
&\qquad\le C_j|\bbx-\bbx'|_{\rm d}^{\g_j}(t-s)^{-\g_j}\big(g_{\lambda}^{(j)}(s,\bbx;t,\bby)+g_{\lambda}^{(j)}(s,\bbx';t,\bby)\big).
\end{split}
\end{align}
\end{enumerate}
All the constants $C_j,\lambda$ appearing above only depend on $\Theta$ and may be different at different places.
\et

Fix $T>0$ and $f\in C^{\infty}_b(\mR^{1+2d})$. Consider the following nonhomogeneous equation:
\begin{align}\label{AG1}
\p_su+\sL^{(0)}_su+f=0,\ \ u(T)=0,
\end{align}
where
$$
\sL^{(0)}_s:=\tr((\sigma\sigma^*)(s,\bbx)\cdot\nabla^2_{x_1})+x_1\cdot\nabla_{x_2}+b_0(s,\bbx)\cdot \nabla_{x_1}.
$$
By the Schauder theory of kinetic equation  (see \cite[Theorem 1]{CHM21}), there is a unique solution $u$ to the above equation with regularities that
\begin{align}\label{Sc1}
u\in C_b([0,T]\times\mR^{2d}),\ \nabla^2_{x_1}u, \nabla_{x_2}u\in \mL^\infty_TC_b(\mR^{2d}).
\end{align}
In particular, by It\^o's formula and \eqref{205}, we have
\begin{align}\label{sl}
u(s,\bbx)=-\int^T_s \bE \left((\p_r u+\sL^{(0)}_ru)(r,\bbX_{r,s}(\bbx))\right)\dif r=\int^T_s P^{(0)}_{r,s}f(r,\bbx)\dif r=:\sI^T_f(s,\bbx).
\end{align}

Next we show some basic estimates about $\sI^T_f$ by Theorem \ref{Th31}.
\bl
Under {\bf (H$^\sigma_b$)}, for $j=0,1$, there exist constants $\lambda\in(0,1), C_j\geq 1$ only depending on $\Theta$ such that for any $t\in[0,1]$,
\begin{align}\label{S6}
\|\nabla^j_{x_1}\sI^t_f\|_{\mL^\infty_t}\leq {C_j} K^{(j)}_{\lambda}(f;t),
\end{align}
and for $\gamma_1\in(0,1)$ and any $s\in[0,t)$, $\delta\in(0,t-s)$ and $\bbx,\bbx'\in\mR^{2d}$,
\begin{align}\label{S4}
|\nabla_{x_1}\sI^t_f(s,\bbx)-\nabla_{x_1}\sI^t_f(s,\bbx')|\lesssim_{C'_1}K^{(1)}_{\lambda}(f;\delta)+|\bbx-\bbx'|_{\rm d}^{\gamma_1}\delta^{-{\gamma_1}}K^{(1)}_{\lambda}(f;1),
\end{align}
where $C'_1=C(\Theta,\gamma_1)$. 
 Moreover, for any $t\in(0,1)$, $N>0$ and each $(s_0,\bbx_0)\in[0,t]\times\mR^{2d}$,
\begin{align}\label{S44}
\lim_{(s,\bbx)\to(s_0,\bbx_0)}\sup_{K^{(1)}_{\lambda}(f; 1)\leq N}|\sI^t_f(s,\bbx)-\sI^t_f(s_0,\bbx_0)|=0.
\end{align}
\el
\begin{proof}
(i) By definition and Theorem \ref{Th31}, we have for $j=0,1$ and $(s,\bbx)\in[0,t]\times\mR^{2d}$,
\begin{align*}
|\nabla^j_{x_1}\sI^t_f(s,\bbx)|
&\le\int^t_s\!\!\!\int_{\mR^{2d}}|\nabla^j_{x_1}p_0|(s,\bbx;r,\bby)|f(r,\bby)|\dif\bby\dif r\\
&\lesssim \int^t_s\!\!\!\int_{\mR^{2d}}g^{(j)}_{\lambda}(s,\bbx;r,\bby)|f(r,\bby)|\dif\bby\dif r\\
&= \int^t_s\!\!\!\int_{\mR^{2d}}\eta^{(j)}_{\lambda}(r-s,\theta_{r,s}(\bbx)-\bby)|f(r,\bby)|\dif\bby\dif r\\
&\leq\sup_{\bbz}\int^{t-s}_0\!\!\!\int_{\mR^{2d}}\eta^{(j)}_{\lambda}(r,\bby)|f(r+s,\bbz-\bby)|\dif\bby\dif r\\
&\leq K^{(j)}_{\lambda}(f; t).
\end{align*}
(ii) Fix $\delta\in(0,t-s)$. By definition we have
$$
|\nabla_{x_1}\sI^t_f(s,\bbx)-\nabla_{x_1}\sI^t_f(s,\bbx')|\leq I_1+I_2,
$$
where
\begin{align*}
I_1&:=\int^t_{s+\delta}\!\int_{\mR^{2d}}|\nabla_{x_1}p_0(s,\bbx;r,\bby)-\nabla_{x_1}p_0(s,\bbx';r,\bby)|\,|f(r,\bby)|\dif\bby\dif r,\\
I_2&:=\int^{s+\delta}_s\!\!\!\int_{\mR^{2d}}|\nabla_{x_1}p_0(s,\bbx;r,\bby)-\nabla_{x_1}p_0(s,\bbx';r,\bby)|\,|f(r,\bby)|\dif\bby\dif r.
\end{align*}
For $I_1$, by \eqref{A1} and the change of variable, we have
\begin{align*}
I_1&\lesssim |\bbx-\bbx'|_{\rm d}^{\gamma_1}\int^t_{s+\delta}\!\int_{\mR^{2d}}(r-s)^{-{\gamma_1}}\Big(g_{\lambda}^{(1)}(s,\bbx;r,\bby)+g_{\lambda}^{(1)}(s,\bbx';r,\bby)\Big)|f(r,\bby)|\dif\bby\dif r\\
&\leq \frac{|\bbx-\bbx'|_{\rm d}^{\gamma_1}}{\delta^{\gamma_1}}\int^t_{s+\delta}\!\int_{\mR^{2d}}\Big(\eta_{\lambda}^{(1)}(r-s,\theta_{r,s}(\bbx)-\bby)+\eta_{\lambda}^{(1)}(r-s,\theta_{r,s}(\bbx')-\bby)\Big)|f(r,\bby)|\dif\bby\dif r\\
&\lesssim \frac{|\bbx-\bbx'|_{\rm d}^{\gamma_1}}{\delta^{\gamma_1}}\sup_{\bbz}\int^t_0\!\!\int_{\mR^{2d}}\eta_{\lambda}^{(1)}(r,\bby)|f(r+s,\bbz-\bby)|\dif\bby\dif r
\leq\frac{|\bbx-\bbx'|_{\rm d}^{\gamma_1}}{\delta^{\gamma_1}}K_{\lambda}^{(1)}(f; t).
\end{align*}
For $I_2$, we similarly have
\begin{align*}
I_2&\lesssim\int^{s+\delta}_s\!\!\!\int_{\mR^{2d}}\Big(g_{\lambda}^{(1)}(s,\bbx;r,\bby)+g_{\lambda}^{(1)}(s,\bbx';r,\bby)\Big)|f(r,\bby)|\dif\bby\dif r\\
&\leq \int^{s+\delta}_s\!\!\!\int_{\mR^{2d}}\Big(\eta_{\lambda}^{(1)}(r-s,\theta_{r,s}(\bbx)-\bby)+\eta_{\lambda}^{(1)}(r-s,\theta_{r,s}(\bbx')-\bby)\Big)|f(r,\bby)|\dif\bby\dif r\\
&\lesssim\sup_{\bbz }\int^\delta_0\!\!\int_{\mR^{2d}}\eta_{\lambda}^{(1)}(r,\bby)|f(r+s,\bbz-\bby)|\dif\bby\dif r\leq K_{\lambda}^{(1)}(f;\delta).
\end{align*}
Hence,
$$
|\nabla_{x_1}\sI^t_f(s,\bbx)-\nabla_{x_1}\sI^t_f(s,\bbx')|\lesssim K_{\lambda}^{(1)}(f;\delta)+|\bbx-\bbx'|_{\rm d}^{\gamma_1}\delta^{-{\gamma_1}}K_{\lambda}^{(1)}(f;t).
$$
(iii) To show \eqref{S44}, by \eqref{S6}, it suffices to show that for each $\bbx\in\mR^{2d}$,
\begin{align}\label{GS1}
\lim_{s\to s_0}\sup_{K^{(1)}_{\lambda}(f; 1)\leq N}|\sI^t_f(s,\bbx)-\sI^t_f(s_0,\bbx)|=0.
\end{align}
Without loss of generality, we assume $0\leq s_0<s<t$. For $\delta\in(0,t-s)$,
by definition we have
\begin{align*}
\sI^t_f(s,\bbx)-\sI^t_f(s_0,\bbx)&=\int^t_{s+\delta}(P^{(0)}_{r,s}-P^{(0)}_{r,s_0})f(r,\bbx)\dif r
+\int^{s+\delta}_{s}P^{(0)}_{r,s}f(r,\bbx)\dif r-\int^{s+\delta}_{s_0}P^{(0)}_{r,s_0}f(r,\bbx)\dif r\\
&=:I_1(s,s_0)+I_2(s,s_0)+I_3(s,s_0).
\end{align*}
For the first term,  if we write $G_s(r,\bbx):=P^{(0)}_{r,s}f(r,\bbx),$ then
\begin{align*}
|I_1(s,s_0)|\leq\int^t_{s+\delta}|G_s(r,\bbx)-P^{(0)}_{s,s_0}G_s(r,\bbx)|\dif r.
\end{align*}
Noting that by \eqref{A1},
$$
|G_s(r,\bbx)-G_s(r,\bbx')|\lesssim |\bbx-\bbx'|_{\rm d}^{\g_0}(r-s)^{-\g_0}\int_{\mR^{2d}}\big(g_{\lambda}^{(0)}(s,\bbx;r,\bby)+g_{\lambda}^{(0)}(s,\bbx';r,\bby)\big)f(r,\bby)\dif \bby,
$$
by definition \eqref{205}, we have
\begin{align*}
&|I_1(s,s_0)|\leq \bE\left(\int^t_{s+\delta}|G_s(r,\bbx)-G_s(r,\bX_{s,s_0}(\bbx))|\dif r\right)\\
&\leq C_\delta\bE\left(|\bbx-\bX_{s,s_0}(\bbx)|_{\rm d}^{\g_0}\int^t_{s+\delta}\int_{\mR^{2d}}\big(g_{\lambda}^{(0)}(s,\bbx;r,\bby)+g_{\lambda}^{(0)}(s,\bX_{s,s_0}(\bbx);r,\bby)\big)f(r,\bby)\dif \bby\right)\\
&\leq C_\delta \bE|\bbx-\bX_{s,s_0}(\bbx)|_{\rm d}^{\g_0}K_{\lambda}^{(0)}(f;t)\leq C_{\delta,\bbx}(s-s_0)^{\gamma_0/2}K_{\lambda}^{(0)}(f;t),
\end{align*}
where the last step is due to \eqref{SDE0} and standard stochastic calculation.
Moreover, by \eqref{6},
$$
|I_2(s,s_0)|\lesssim K_{\lambda}^{(0)}(f;\delta),\ \ |I_3(s,s_0)|\lesssim K_{\lambda}^{(0)}(f;s-s_0+\delta).
$$
Hence, by \eqref{AA1} with $\beta'=0$ and $\beta=1$,
$$
\sup_{K^{(1)}_{\lambda}(f; 1)\leq N}|\sI^t_f(s,\bbx)-\sI^t_f(s_0,\bbx)|\leq C_{N,\delta,\bbx} (s-s_0)^{\gamma_0/2}
+C_N(\delta^{1/2}+(s-s_0+\delta)^{1/2}),
$$
which in turn yields \eqref{GS1}. 
The proof is complete.
\end{proof}
We need the following interpolation limit result (cf. \cite[Chapter 3]{Ky96}).
\bl\label{Le33}
Let $E$ be an index set and $D\subset\mR^m$ a bounded open set.
Let $f_n(t,x): E\times D\to\mR$ be a sequence of differentiable functions in $x$. Assume that for some $f: E\times D\to\mR$,
\begin{align}\label{AW1}
\lim_{n\to\infty}\sup_{(t,x)\in E\times D}|f_n(t,x)-f(t,x)|=0,
\end{align}
and
\begin{align}\label{AW2}
\sup_{(t,x)\in E\times D}\sup_n|\nabla f_n(t,x)|<\infty, \ \ \lim_{|x-y|\to 0}\sup_{t\in E}\sup_n|\nabla f_n(t,x)-\nabla f_n(t,y)|=0.
\end{align}
Then $x\mapsto\nabla f(t,x)$ is uniformly continuous in $t\in E$ and
\begin{align}\label{AW24}
\lim_{n\to\infty}\sup_{(t,x)\in E\times D}|\nabla f_n(t,x)-\nabla f(t,x)|=0.
\end{align}
\el
\begin{proof}
By \eqref{AW2} and Ascoli-Arzela's theorem, there is a subsequence $n_k$ and an $\mR^m$-valued function $H: E\times D\to\mR$, which is continuous in $x$ and uniformly in $t\in E$, so that 
$$
\lim_{k\to\infty}\sup_{(t,x)\in E\times D}|\nabla f_{n_k}(t,x)-H(t,x)|=0,
$$
which together with \eqref{AW1} yields 
$$
H(t,x)=\nabla f(t,x).
$$
Now we prove \eqref{AW24}.
For simplicity of notation, we write $g_n:=f_n-f$. Fix $i=1,\cdots, m$, $x\in D$ and $r<{\rm dis}(x, \p D)$. Let $e_i=(0,\cdots, 1, \cdots,0)$ be the unit vector in $\mR^m$.
By the mean valued theorem, there is a point $y\in B_r(x)$ such that
\begin{align}\label{AW25}
g_n(t,x)-g_n(t,x+re_i)=r\p_{x_i}g_n(t,y).
\end{align}
For any $\eps\in(0,1)$, by the assumption, we can choose $r$ small enough so that
$$
\sup_{|x-y|\leq r}\sup_{t\in E}\sup_n|\p_{x_i} g_n(t,x)-\p_{x_i} g_n(t,y)|\leq \eps,
$$
and by \eqref{AW25},
$$
|\p_{x_i} g_n(t,x)|\leq|\p_{x_i} g_n(t,x)-\p_{x_i} g_n(t,y)|+|\p_{x_i}g_n(t,y)|\leq\eps+2\sup_{z\in B_r(x)}|g_n(t,z)|/r.
$$
From these, we derive the desired limit by letting $n\to\infty$ and \eqref{AW1}.
\end{proof}

\bl\label{Le34}
Let $(f_n)_{n\in\mN}\subset\mK^{\rm o}_0$. Suppose that for any $\lambda>0$,
$$
\lim_{\delta\to 0}\sup_nK^{(0)}_\lambda(f_n; \delta)=0,
$$
and for any $t>0$ and bounded domain $D\subset\mR^{2d}$,
$$
\lim_{n\to\infty}\int^t_0\!\!\int_D|f_n(r,\bby)|\dif \bby\dif r=0.
$$
Then for any $t>0$ and  bounded domain $D\subset\mR^{2d}$,
\begin{align}\label{S5}
\lim_{n\to\infty}\sup_{(s,\bbx)\in[0,t]\times D}|\sI^t_{f_n}(s,\bbx)|=0.
\end{align}
\el
\begin{proof}
Fix $(s,\bbx)\in[0,t)\times D$ and $\delta\in(0,t-s)$.
By definition, we have
$$
|\sI^t_{f_n}(s,\bbx)|\leq \int^{s+\delta}_s|P^{(0)}_{r,s}f_n(r,\bbx)|\dif r+\int^t_{s+\delta}|P^{(0)}_{r,s}f_n(r,\bbx)|\dif r=:I_1(s,\bbx)+I_2(s,\bbx).
$$
For $I_1(s,\bbx)$, by  \eqref{S6} we have
$$
I_1(s,\bbx)\lesssim K^{(0)}_\lambda(f_n;\delta)\leq K^{(0)}_\lambda(f;\delta).
$$
For $I_2(s,\bbx)$, let $R$ be large enough so that $\{\theta_{r,s}(\bbx),r\in[s,t]\}\subset B_R$.  Then by \eqref{6}, we have
\begin{align*}
I_2(s,\bbx)&\lesssim \int^t_{s+\delta}\int_{\mR^{2d}}\eta^{(0)}_\lambda(r-s,\theta_{r,s}(\bbx)-\bby)|f_n(r,\bby)|\dif \bby\dif r\\
&\leq \sup_{\bbz\in B_R}\int^t_{s+\delta}\int_{\mR^{2d}}\eta^{(0)}_\lambda(r-s,\bbz-\bby)|f_n(r,\bby)|\dif \bby\dif r\\
&\leq \sup_{\bbz\in B_R}\int^t_{s+\delta}\int_{B_{2R}}\eta^{(0)}_\lambda(r-s,\bbz-\bby)|f_n(r,\bby)|\dif \bby\dif r\\
&\quad+\sup_{\bbz\in B_R}\int^t_{s+\delta}\int_{B^c_{2R}}\eta^{(0)}_\lambda(r-s,\bbz-\bby)|f_n(r,\bby)|\dif \bby\dif r\\
&=:I_{21}(s,\bbx)+I_{22}(s,\bbx).
\end{align*}
For $I_{21}(s,\bbx)$, since $\eta^{(0)}_\lambda(r-s,\bbz-\bby)\le(r-s)^{-2d}$, we have
$$
I_{21}(s,\bbx)\leq \delta^{-2d}\int^t_{s+\delta}\int_{B_{2R}}|f_n(r,\bby)|\dif \bby\dif r.
$$
For $I_{22}(s,\bbx)$, by the definition of $\eta^{(0)}_\lambda$ we have
\begin{align*}
I_{22}(s,\bbx)&\leq
\sup_{\bbz\in B_R}\int^t_{s+\delta}\int_{B^c_{2R}}\eta^{(0)}_{\lambda/2}(r-s,\bbz-\bby)|f_n(r,\bby)|\dif \bby\dif r\\
&\quad\times \sup_{r\in[s+\delta,t]}\sup_{\bbz\in B_R}\sup_{\bby\in B^c_{2R}}\e^{-\lambda |\mT_{r-s}(\bbz-\bby)|^2/2}\\
&\leq  K^{(0)}_{\lambda/2}(f_n; t)\e^{-C_\delta R^2}\leq  K^{(0)}_{\lambda/2}(f; t)\e^{-C_\delta R^2},
\end{align*}
where $C_\delta$ is independent of $R,n$. Combining the above estimates, we obtain 
$$
|\sI^t_{f_n}(s,\bbx)|\leq K^{(0)}_\lambda(f;\delta)+\delta^{-2d}\int^t_{s+\delta}\int_{B_{2R}}|f_n(r,\bby)|\dif \bby\dif r+K^{(0)}_{\lambda/2}(f; t)\e^{-C_\delta R^2}.
$$
Firstly letting $n\to\infty$, then $R\to\infty$ and finally $\delta\to 0$, we get the desired limit.
\end{proof}

Due to $|b_1|\in\mK^{\rm o}_1$, 
for the $\lambda, C_1$ from \eqref{S6}, we may choose $T\in(0,1)$ small enough so that
\begin{align}\label{13}
K^{(1)}_{\lambda}(|b_1|;T)\leq 1/(2C_1).
\end{align}
In the following, we shall fix such a time $T$ and consider the integral equation:
\begin{align}\label{34}
u(s,\bbx)=\sI^T_f(s,\bbx)+\sI^T_{b_1\cdot\nabla_{x_1}u}(s,\bbx),\ \ s\in[0,T].
\end{align}
We now establish the following main result of this section.
\bt\label{37}
Under {\bf (H$^\sigma_b$)}, for any $f\in\mK^{\rm o}_1$, there exists a  unique solution $u\in C_b([0,T]\times\mR^{2d})$ to integral equation \eqref{34} with regularity estimate
\begin{align}\label{14}
\|u\|_{\mL^\infty_{T}}+\|\nabla_{x_1}u\|_{\mL^\infty_{T}}\le C_3 K^{(1)}_{\lambda}(f; T),
\end{align}
where $C_3=C_3(\Theta)>0$ and $\lambda=\lambda(\Theta)\in(0,1)$. 
Moreover, for each $t,\bbx$, we have
\begin{align}\label{35}
\lim_{n\to\infty}u_n(t,\bbx)=u(t,\bbx),
\end{align}
where  $u_n$ is the solution of \eqref{34} corresponding to $(b^n_1, f_n)=(b_1*\phi_n, f*\phi_n)$.
\et
\begin{proof}
(i) By the standard Picard iteration, it suffices to show the estimate \eqref{14}.
For $j=0,1$, by \eqref{S6} we have
\begin{align}\label{nn1}
\|\nabla^j_{x_1}u\|_{\mL^\infty_T}&\lesssim_{C_j} K^{(j)}_{\lambda}(b_1\cdot\nabla_{x_1}u; T)+K^{(j)}_{\lambda}(f;T)\no\\
&\lesssim_{C_j} K^{(j)}_{\lambda}(|b_1|; T)\|\nabla_{x_1}u\|_{\mL^\infty_T}+K^{(j)}_{\lambda}(f;T).
\end{align}
For $j=1$, using \eqref{13} we obtain
$$
\|\nabla_{x_1}u\|_{\mL^\infty_{T}}\leq 2C_1K^{(1)}_{\lambda}(f;T).
$$
Substituting this back into \eqref{nn1} with $j=0$ and by \eqref{AA1}, we get
$$
\|u\|_{\mL^\infty_{T}}\leq 2C_0C_1K^{(1)}_{\lambda}(f;T)+C_0K^{(0)}_{\lambda}(f;T)\leq C_3K^{(1)}_{\lambda}(f;T).
$$
(ii)  Now we prove \eqref{35}. Let $u_n$ solve the following integral equation:
\begin{align}\label{304}
u_n(s,\bbx)=\sI^T_{f_n}(s,\bbx)+\sI^T_{b^n_1\cdot\nabla_{x_1}u_n}(s,\bbx),\ \ s\in[0,T].
\end{align}
By \eqref{14}, we have
$$
\sup_{n\in\mN}\Big(\|u_n\|_{\mL^\infty_{T}}+\|\nabla_{x_1}u_n\|_{\mL^\infty_{T}}\Big)<\infty,
$$
and by \eqref{S44}, 
$$
\lim_{(s,\bbx)\to(s_0,\bbx_0)}\sup_n|u_n(s,\bbx)-u_n(s_0,\bbx_0)|=0.
$$
Therefore, by Ascoli-Arzela's theorem, there is a subsequence $\{n_k\}$ and $u\in C_b([0,t]\times\mR^{2d})$ such that for each bounded domain $D\subset\mR^{2d}$,
$$
\lim_{k\to\infty}\sup_{s\in[0,t],\bbx\in D}|u_{n_k}(s,\bbx)-u(s,\bbx)|=0.
$$
It remains to show that $u$ satisfies \eqref{34}. By \eqref{S4}, for any $\delta\in(0,t-s)$, we have
\begin{align}\label{1000}
|\nabla_{x_1}u_n(s,\bbx)-\nabla_{x_1}u_n(s,\bbx')|
&\lesssim|\bbx-\bbx'|_{\bf d}^{\gamma_1}\delta^{-{\gamma_1}}K^{(1)}_{\lambda}(f_n+b^n_1\cdot\nabla_{x_1}u_n;1)
+K^{(1)}_{\lambda}(f_n+b^n_1\cdot\nabla_{x_1}u_n;\delta)\nonumber\\
&\lesssim|\bbx-\bbx'|_{\bf d}^{\gamma_1}\delta^{-{\gamma_1}}
\Big(K^{(1)}_{\lambda}(|f_n|;1)+\|\nabla_{x_1}u_n\|_\infty K^{(1)}_{\lambda}(|b^n_1|;1)\Big)\nonumber\\
&\quad+K^{(1)}_{\lambda}(|f_n|;\delta)+\|\nabla_{x_1}u_n\|_\infty K^{(1)}_{\lambda}(|b^n_1|;\delta)\nonumber\\
&\lesssim|\bbx-\bbx'|_{\bf d}^{\gamma_1}\delta^{-{\gamma_1}}+K^{(1)}_{\lambda}(|f|;\delta)+ K^{(1)}_{\lambda}(|b_1|;\delta),
\end{align}
where the implicit constant does not depend on $n,\delta$ and $\bbx,\bbx'$. Thus, by Lemma \ref{Le33}, 
$$
\lim_{k\to\infty}\sup_{s\in[0,t],\bbx\in D}|\nabla_{x_1}u_{m_k}(s,\bbx)-\nabla_{x_1}u(s,\bbx)|=0.
$$
Now, by Lemma \ref{Le34} and taking limits for \eqref{304}, one sees that $u$ is a solution of \eqref{34}. Finally, by the uniqueness, the whole sequence converges.
\end{proof}

\section{Well-posedness for SDE with singular drift}

In this section, we aim to demonstrate the well-posedness of SDE with drift $b_1$ belonging to the Kato class and provide a proof for the first part of Theorem \ref{Main1}. Moreover, throughout the following discussion, we still fix $T$ as defined in \eqref{13}.

We first establish the following crucial Krylov estimate based on \eqref{S6}.

\bl\label{Le41}
Let $\bbX$ be a weak solution of SDE \eqref{1} in the sense of Definition \ref{weak}. 
Under {\bf (H$^\sigma_b$)}, there are constants $C_3>0$ and $\lambda\in(0,1)$ depending only on $\Theta$ such that for any $0\leq f\in \mK^{\rm o}_1$,
\begin{align}\label{AA05}
\bE\left(\int_0^tf(s,\bbX_s)\dif s\right)\le  C_3 K^{(1)}_{\lambda}(f;t),\ \ \forall t\in(0,T].
\end{align}
\el
\begin{proof}
Fix $t\in(0,T]$ and the starting point $\bbx\in\mR^{2d}$ of SDE \eqref{25}.
Firstly we consider $0\leq f\in C^\infty_b(\mR^{1+2d})$ and let $u(s,\cdot):=\sI^t_f(s,\cdot)$.
By \eqref{Sc1}, It\^o's formula and \eqref{AG1}, we have for any $r\in[0,t]$,
\begin{align*}
u(r,\bbX_r)&=u(0,\bbx)+\int_0^r(\p_s u+\sL_s u)(s,\bbX_s)\dif s+\sqrt{2}\int_0^r(\sigma^*\cdot\nabla_{x_1}u)(s,\bbX_s)\dif W_s\\
&=u(0,\bbx)+\int_0^r(b_1\cdot\nabla_{x_1}u-f)(s,\bbX_s)\dif s+\sqrt{2}\int_0^r(\sigma^*\cdot\nabla_{x_1}u)(s,\bbX_s)\dif W_s.
\end{align*}
Let $R>|\bbx|$ and define a stopping time
$$
\tau_R:=\inf\bigg\{r>0: \int^r_0|b_1(s,\bbX_s)|\dif s\geq R\mbox{ or } |\bbX_r|\geq R\bigg\}.
$$
By BDG's inequality, it is easy to see that for some $C_2=C_2(\Theta)>0$,
$$
\bE\left(\int_0^{t\wedge\tau_R} f(s,\bbX_s)\dif s\right)\leq 2\|u\|_{\mL^\infty_t}+\|\nabla_{x_1}u\|_{\mL^\infty_t}\left[\bE\left(\int_0^{t\wedge\tau_R}|b_1(s,\bbX_s)|\dif s\right)+C_2\right].
$$
By \eqref{S6}, we further have
\begin{align}\label{AA03}
\bE\left(\int_0^{t\wedge\tau_R} f(s,\bbX_s)\dif s\right)&\leq 2C_0K^{(0)}_{\lambda}(f;t)
+C_1K^{(1)}_{\lambda}(f;t)\left[\bE\left(\int_0^{t\wedge\tau_R}|b_1(s,\bbX_s)|\dif s\right)+C_2\right].
\end{align}
Let $b^n_1:=b_1*\phi_n\in C^\infty_b(\mR^{1+2d})$. Then for Lebesgue almost all $(s,\bbx)$,
$$
\lim_{n\to\infty}b^n_1(s,\bbx)=b_1(s,\bbx).
$$
Since the law of $(s,\bbX_s)$ is absolutely continuous with respect to the Lebesgue measure, we have
$$
\lim_{n\to\infty}b^n_1(s,\bbX_s)=b_1(s,\bbX_s),\  \ a.s.
$$
Hence, by Fatou's lemma and \eqref{AA03} with $f=b^n_1$, 
\begin{align*}
&\bE\left(\int_0^{t\wedge\tau_R} |b_1(s,\bbX_s)|\dif s\right)\le\varliminf_{n\to\infty}\bE\left(\int_0^{t\wedge\tau_R} |b^n_1(s,\bbX_s)|\dif s\right)\\
&\quad\leq \sup_n2C_0K^{(0)}_{\lambda}(|b^n_1|;t)+ \sup_nC_1K^{(1)}_{\lambda}(|b^n_1|;t)\left[\bE\left(\int_0^{t\wedge\tau_R}|b_1(s,\bbX_s)|\dif s\right)+C_2\right]\\
&\quad\leq 2C_0K^{(0)}_{\lambda}(|b_1|;t)+ \frac12\left[\bE\left(\int_0^{t\wedge\tau_R}|b_1(s,\bbX_s)|\dif s\right)+C_2\right],
\end{align*}
where the last step is due to  \eqref{CC1} and \eqref{13}.
This implies that for any $R>0$,
\begin{align}\label{32}
\bE\left(\int_0^{t\wedge\tau_R} |b_1(s,\bbX_s)|\dif s\right)\leq 4C_0K^{(0)}_{\lambda}(|b_1|;t)+C_2.
\end{align}
Substituting this into \eqref{AA03} and letting $R\to\infty$, we obtain the desired estimate for $f\in C^\infty_b(\mR^{1+2d})$.

Next, we consider the general case $0\leq f\in \mK^{\rm o}_1$. Let 
$$
f_n:=f*\phi_n\in C^\infty_b(\mR^{1+2d}).
$$ 
By Fatou's lemma and what we have proven, we have
\begin{align*}
\bE\left(\int_0^t f(s,\bbX_s)\dif s\right)
\le\varliminf_{n\to\infty}\bE\left(\int_0^t f_n(s,\bbX_s)\dif s\right)
\le  C_3 \sup_nK^{(1)}_{\lambda}(f_n;t)\stackrel{\eqref{CC1}}{\le} C_3 K^{(1)}_{\lambda}(f;t).
\end{align*}
The proof is complete.
\end{proof}

Next we consider the approximating SDE with smooth coefficients:
\begin{align}\label{SDEN}
\dif\bbX^{n}_t=(B_0+B_1^n)(t,\bbX^n_t)\dif t+\Sigma_n\(t,\bbX_t^n\)\dif W_t,\quad \bbX^{n}_0=\bbx,
\end{align}
where
\begin{align}\label{SX3}
B_1^n(t,\bbx)=(b_1*\phi_n(t,\bbx), 0)^*,\ \ \Sigma_n(t,\bbx)=(\sqrt{2}(\sigma*\phi_n)(t,\bbx),0_{d\times d})^*.
\end{align}
Since the coefficients are global Lipschitz, there is a unique strong solution to the above SDE.
Let $\mP_n$ be the law of $\bbX^n$ in the space $\mC_T$. We have
\bt\label{Th43}
Under {\bf (H$^\sigma_b$)}, for any $\bbx\in\mR^{2d}$, there is a unique weak solution $\{\bbX_t(\bbx)\}_{t\in[0,T]}$ to SDE \eqref{25} starting from $\bbx$. Moreover, 
$\mP_n$ weakly converges to the law of $\bbX_\cdot$ in $\mC_T$.
\et
\begin{proof}
({\bf Uniqueness})
Let $f\in C^\infty_b(\mR^{1+2d})$. Consider the following backward PDE with mollifying $b^n_1$:
\begin{align}\label{eq3}
\p_su_n+\sL_s^{(0)}u_n+b_1^n\cdot\nabla_{x_1}u_n+f=0,\ \ u_n(T,\bbx)=0.
\end{align}
By the Schauder theory (see \cite{CHM21}), there is a unique solution $u_n$ with the regularity
$$
\nabla^2_{x_1}u_n,\nabla_{x_2}u_n\in \mL^\infty_TC_b(\mR^{2d}).
$$ 
Let $\bbX_t:=\bbX_t(\bbx)$ be a solution of SDE \eqref{25} starting from $\bbx$.
By It\^o's formula and \eqref{eq3}, we have
$$
u_n(t,\bbX_t)=u_n(0,\bbx)-\int_0^tf(s,\bbX_s)\dif s+\int_0^t((b_1-b_1^{n})\cdot\nabla_{x_1}u_n)(s,\bbX_s)\dif s+\mbox{ a martingale}.
$$
Let $R>|\bbx|$ and define a stopping time
$$
\tau_R:=\inf\big\{t>0: |\bbX_t|\geq R\big\}.
$$
By the optional stopping theorem, we have
\begin{align}\label{AG12}
\bE u_n(T\wedge\tau_R,\bbX_{T\wedge\tau_R})=u_n(0,\bbx)+\bE\left(\int_0^{T\wedge\tau_R}f(s,\bbX_s)\dif s\right)+J^{R}_n,
\end{align}
where
$$
J^{R}_n:=\bE\left(\int_0^{T\wedge\tau_R}((b_1-b_1^{n})\cdot\nabla_{x_1}u_n)(s,\bbX_s)\dif s\right).
$$
For $J^{R}_n$, by \eqref{CC1} and \eqref{AA05}, we have 
\begin{align*}
|J_n^R|&\le\|\nabla_{x_1}u_n\|_{\mL_T^\infty}\bE\left(\int_0^{T\wedge\tau_R}|b_1-b_1^{n}|(s,\bbX_s)\dif s\right)\\
&\le\|\nabla_{x_1}u_n\|_{\mL_T^\infty}\bE\left(\int_0^{T}\1_{|\bbX_s|\leq R}|b_1-b_1^{n}|(s,\bbX_s)\dif s\right)\\
&\lesssim K^{(1)}_{\lambda}(\1_{[0,T]\times B_R}|b_1-b^n_1|; T),
\end{align*}
which implies by Lemma \ref{Le24} that
$$
\lim_{n\to\infty}|J_n^R|= 0.
$$
By Duhamel's formula, $u_n$ also satisfies the following integral equation:
$$
u_n(s,\bbx)=\sI^T_f(s,\bbx)+\sI^T_{b^n_1\cdot\nabla_{x_1}u_n}(s,\bbx), 
$$
where $\sI^T_f$ is defined in \eqref{sl}.
Let $u\in C_b([0,T]\times\mR^{2d})$ be the unique  solution of integral equation \eqref{34}. 
By Theorem \ref{37}, we have 
for each $t,\bbx$,
\begin{align}\label{AG2}
\lim_{n\to\infty}u_n(t,\bbx)=u(t,\bbx).
\end{align}
Thus by the dominated convergence theorem and taking limits $n\to\infty$ for \eqref{AG12}, we get
\begin{align*}
\bE u(T\wedge\tau_R,\bbX_{T\wedge\tau_R})=u(0,\bbx)-\bE\left(\int_0^{T\wedge\tau_R}f(s,\bbX_s)\dif s\right).
\end{align*}
Finally, letting $R\to\infty$ and noting $u(T,\bbx)=0$, we obtain
$$
u(0,\bbx)=\bE\left(\int_0^{T}f(s,\bbX_s)\dif s\right).
$$
From this, it is now standard to derive  the uniqueness  (see \cite[Theorem 4.4.3]{EK86}).

({\bf Existence}) We use the standard weak convergence method. Note that
$$
\bbX^n_t=\bbx+\int_0^tB_0(s,\bbX^n_s)\dif s+\int_0^tB^n_1(s,\bbX^n_s)\dif s+\int_0^t\Sigma_n(s,\bbX_s^n)\dif W_s.
$$
We first show the tightness of $U^n_t:=\int^t_0 B^n_1(s,\bbX^n_s)\dif s=\Big(\int^t_0 b^n_1(s,\bbX^n_s)\dif s, 0\Big)^*, n\in\mN$ in $\mC_T$.
For any stopping time $\tau\le T$, by the strong Markov property and \eqref{AA05}, we have
\begin{align*}
\bE|U^n_{\tau+\delta}-U^n_\tau|
&\leq\bE\left(\int_\tau^{\tau+\delta}|b_1^n(s,\bbX_{s}^n(\bbx))|\dif s\right)\\
&\le\sup_{t,\bby}\bE\left(\int_t^{t+\delta}|b_1^n(s,\bbX_s^n(\bby))|\dif s\right)\\
&\leq C_3 K^{(1)}_{\lambda}(|b_1^n|;\delta)\leq C_3 K^{(1)}_{\lambda}(|b_1|;\delta).
\end{align*}
By \cite[Lemma 2.7]{ZZ18}, we have
\begin{align*}
\bE\left(\sup_{t\in[0,T]}|U^n_{t+\delta}-U^n_t|^{1/2}\right)\lesssim(K^{(1)}_{\lambda}(|b_1|;\delta)) ^{1/2},
\end{align*}
which implies that the law of $U^n_\cdot$ in $\mC_T$ is tight.

Moreover, we have
\begin{align*}
|\bbX_t^n|
\leq|\bbx|+\int_0^t|B_0(s,\bbX^n_s)|\dif s+\int_0^t|b^n_1(s,\bbX^n_s)|\dif s+\left|\int_0^t\Sigma_n(s,\bbX_s^n)\dif W_s\right|.
\end{align*}
Since $B_0$ is linear growth and $\Sigma_n$ is bounded, by BDG's inequality, we have
\begin{align*}
\bE\left(\sup_{s\in[0,t]}|\bbX_s^n|\right)
&\leq|\bbx|+C\int_0^t(1+\bE|\bbX^n_s|)\dif s+\bE\left(\int_0^t|b^n_1(s,\bbX^n_s)|\dif s\right)+C\sqrt t,
\end{align*}
which implies by Gronwall's inequality and \eqref{AA05} that
\begin{align}\label{Lim9}
\bE\left(\sup_{s\in[0,t]}|\bbX_s^n|\right)
\lesssim 1+|\bbx|+\bE\left(\int_0^t|b^n_1(s,\bbX^n_s)|\dif s\right)\lesssim 1+|\bbx|.
\end{align}
Thus,
\begin{align*}
\bE\left(\sup_{t\in[0,T]}\int^{t+\delta}_t |B_0(s,\bbX^n_s)|\dif s\right)\leq C\delta \bE\left(1+\sup_{s\in[0,t]}|\bbX_s^n|\right)\leq C(1+|\bbx|)\delta.
\end{align*}
Consequently, the law of $(\int^\cdot_0 B_0(s,\bbX^n_s)\dif s)_{n\in\mN}$ in $\mC_T$ is tight.
Moreover, by BDG's inequality, it is also easy to see that the law of $(\int^\cdot_0 \Sigma_n(s,\bbX^n_s)\dif W_s)_{n\in\mN}$ in $\mC_T$ is tight.
Combining these tightness, we obtain the tightness of $(\mQ^n)_{n\in\mN}:=(\bP\circ(\bbX^n, W)^{-1})_{n\in\mN}$ in $\mC_T\times\mC_T$.
Hence, there is a subsequence still denoted by $n$ so that $\mQ^n$ weakly converges to some probability measure $\mQ$. By Skorokhod's representation theorem, there are probablity space $(\tilde{\Omega},\tilde{\sF},\tilde{\mP})$ and $\mC_T\times\mC_T$-valued random variables $(\tilde{\bbX}^n,\tilde{W}^n)$ and  $(\tilde{\bbX},\tilde{W})$ defined on it such that 
\begin{align}\label{40}
\tilde{\mP}\circ(\tilde{\bbX}^n,\tilde{W}^n)^{-1}=\mQ^n,\ \ \tilde{\mP}\circ(\tilde{\bbX},\tilde{W})^{-1}=\mQ,
\end{align}
\begin{align}\label{410}
(\tilde{\bbX}^n,\tilde{W}^n)\to(\tilde{\bbX},\tilde{W}),\ \ \tilde{\mP}-a.s.
\end{align}
and 
\begin{align}\label{SDE1}
\dif\tilde\bbX^{n}_t=(B_0+B_1^n)(t,\tilde{\bbX}^n_t)\dif t+\Sigma_n\(t,\tilde{\bbX}_t^n\)\dif \tilde W^n_t.
\end{align}
To prove that $\tilde\bbX_\cdot$ is a weak solution, it suffices to show  that for each $t\in[0,T]$,
\begin{align}\label{eq22}
\int_0^t\Sigma_n(s,\tilde{\bbX}^n_s)\dif \tilde W^n_s\stackrel{n\to\infty}{\to} \int^t_0\Sigma(s,\tilde{\bbX}_s)\dif \tilde W_s\mbox{ in probability},
\end{align}
and
\begin{align}\label{eq2}
\int_0^tb_1^n(s,\tilde{\bbX}^n_s)\dif s\stackrel{n\to\infty}{\to} \int^t_0b_1(s,\tilde{\bbX}_s)\dif s\mbox{ in probability}.
\end{align}
For limit \eqref{eq22}, it follows by \eqref{410} and \cite[p.32]{SK82}, we omit the details.
Next we look at \eqref{eq2}.
For fixed $m$, by the dominated convergence theorem, we have
$$
\lim_{n\to\infty}\tilde\mE\left(\int_0^T |b^m_1(s,\tilde{\bbX}^n_s)-b^m_1(s,\tilde{\bbX}_s)|\dif s\right)=0.
$$
For $R>|\bbx|$ and $n\in\mN$, define a stopping time
$$
\tau^n_R:=\inf\big\{t>0: |\tilde\bbX^n_t|\geq R\big\}.
$$
By \eqref{AA05}, we have
$$
\tilde\mE\left(\int_0^{T\wedge\tau^n_R} |b_1^m-b_1|(s,\tilde{\bbX}^n_s)\dif s\right)
\leq  C_3K^{(1)}_{\lambda}(\1_{[0,T]\times B_R}|b_1^m-b_1|;T),
$$
which implies by Lemma \ref{Le24} that for fixed $R>|\bbx|$,
$$
\lim_{m\to\infty}\sup_n\tilde\mE\left(\int_0^{T\wedge\tau^n_R} |b_1^m-b_1|(s,\tilde{\bbX}^n_s)\dif s\right)=0.
$$
On the other hand, by \eqref{Lim9} we also have
$$
\lim_{R\to\infty}\sup_n\tilde \mP(\tau^n_R>T)
=\lim_{R\to\infty}\sup_n\tilde \mP\left(\sup_{t\in[0,T]}|\tilde\bbX^n_t|\geq R\right)
\leq\lim_{R\to\infty}\sup_n\bE\left(\sup_{t\in[0,T]}|\bbX^n_t|\right)/R=0.
$$
Combining the above limits, we obtain \eqref{eq2}.
The proof is complete.
\end{proof}
\br
Although the existence and uniqueness of a solution is only obtained on a small time interval, 
it can be extended to any large time by a standard patching up technique (cf. \cite{SV97}).  
\er

\section{Two-sided estimates of heat kernel}

In this section we devote to proving the existence of density and two-sided estimates \eqref{AX1}.
Our strategy is as follows: We first prove \eqref{AX1} for the smooth approximation SDE \eqref{SDEN}. Then by taking weak limits, we immediately obtain the desired estimates.
More precisely, consider the approximation SDE:
\begin{align}\label{ASDE}
\dif\bbX^{n}_{t,s}=B_n(t,\bbX^n_{t,s})\dif t+\Sigma_n\(t,\bbX_{t,s}^n\)\dif W_t,\ \ t\geq s,
\end{align}
where $B_n=B_0+B^n_1$ and $\Sigma_n$ are given in \eqref{SX3}.
By \cite{DM10} or \cite[Theorem 1.1]{CMPZ23} , $\bbX^n_{t,s}$ admits a  density $p_n(s,x;t,y)$ that enjoys the following two-sided estimates:
$$
C_0(t-s)^{-2d}\e^{-\lambda_0|\mT_{t-s}(\t_{t,s}(\bbx)-\bby)|^2}\le p_n(s,\bbx;t,\bby) \le C_1(t-s)^{-2d}\e^{-\lambda_1|\mT_{t-s}(\t_{t,s}(\bbx)-\bby)|^2},
$$
where, in \cite{DM10} and \cite{CMPZ23}, the constants $C_0,C_1, \lambda_0, \lambda_1$ depend on $\kappa_0,\gamma_0$ in {\bf (H$^\sigma_b$)} and the Lipschitz or H\"older constants of $B_n$.

The aim below is to show that the above constants only depend on $\Theta$ in \eqref{Para} and the quantity $K^{(1)}_{\lambda}(|b_1|^\gamma,1)$, where $\gamma>1$ and for $i=0,1$, $\lambda_i=\lambda_i(\Theta)>0$.
Once we show this, then noting that the above two-sided estimates are equivalent to that for any $f\in C_b(\mR^{2d})$,
\begin{align*}
&C_0(t-s)^{-2d}\int_{\mR^{2d}}\e^{-\lambda_0|\mT_{t-s}(\t_{t,s}(\bbx)-\bby)|^2}f(\bby)\dif\bby\\
&\qquad\qquad\leq
\bE f(\bbX^n_{t,s}(\bbx))\leq C_1(t-s)^{-2d}\int_{\mR^{2d}}\e^{-\lambda_1|\mT_{t-s}(\t_{t,s}(\bbx)-\bby)|^2}f(\bby)\dif\bby,
\end{align*}
by Theorem \ref{Th43} and taking weak limits $n\to\infty$, we immediately obtain \eqref{AX1}.
In the following we always assume {\bf (H$^\sigma_b$)} and 
\begin{align}\label{AD1}
b_1,\sigma\in \cap_{T>0}\mL^\infty_TC^\infty_b(\mR^{2d}).
\end{align}
We will use the convention that all the constants appearing below only depend on $\Theta$ in  \eqref{Para} and the quantity $K^{(1)}_{\lambda}(|b_1|^\gamma,1)$.
We shall fix $T>0$ as in \eqref{13}.

Let $p_0$ be the heat kernel of operator $\sL^{(0)}$ with regular drift $b_0$ and $p$ be the heat kernel of operator $\sL$.
By the Duhamel representation (see \cite[Subsection 3.1]{CMPZ23}), we have
\begin{align}\label{3}
p(s,\bbx;t,\bby)=p_0(s,\bbx;t,\bby)+(p\otimes H)(s,\bbx;t,\bby),
\end{align}
where 
$$
p\otimes H(s,\bbx;t,\bby)=\int_s^t\!\!\!\int_{\mR^{2d}}p(s,\bbx;r,\bbz)H(r,\bbz;t,\bby)\dif\bbz\dif r,
$$
and
$$
H(s,\bbx;t,\bby)=b_1(s,\bbx)\cdot\nabla_{x_1}p_0(s,\bbx;t,\bby).
$$
For $N\ge 2$, iterating equation \eqref{3} $N-1$ times, we obtain
\begin{align}\label{4}
p(s,\bbx;t,\bby)=p_0(s,\bbx;t,\bby)+\sum_{j=1}^{N-1}(p_0\otimes H^{\otimes j})(s,\bbx;t,\bby)+(p\otimes H^{\otimes N})(s,\bbx;t,\bby).
\end{align}
For simplicity of notations, we write
\begin{align}\label{40}
\Lambda_\lambda(t):=K^{(1)}_\lambda(|b_1|; t).
\end{align}

The following lemma provides an estimate for the summation.
\bp
For any $T\in(0,1]$ and $N\in\mN$, there are constants $C_N>0,\l_N\in (0,1)$ depending only on $\Theta$ and $N$ such that for all $(s,\bbx;t,\bby)\in\mD_T$,
\begin{align}
\sum_{j=1}^{N}|p_0\otimes H^{\otimes j}|(s,\bbx;t,\bby)\le C_N\Lambda^N_{\lambda_N}(t-s) g_{\l_N}(s,\bbx;t,\bby),
\end{align}
where $\l_N\downarrow 0$ as $N\to\infty$.
\ep
\begin{proof}
First of all, we use induction to show that
\begin{align}\label{5h}
|H^{\otimes N}(s,\bbx;t,\bby)|\lesssim  \Lambda^{N-1}_{\lambda_{N-1}}(t-s)|b_1(s,\bbx)| g^{(1)}_{\l_N}(s,\bbx;t,\bby),
\end{align}
where $\l_N=\k^{N-1}\lambda$, and $\lambda,\k\in(0,1)$ are from \eqref{5} and  \eqref{7}. 
For $N=1$, 
by \eqref{5} we have
$$
|H|(s,\bbx;t,\bby)
\lesssim\left|b_1(s,\bbx)\right|g^{(1)}_{\lambda}(s,\bbx;t,\bby).
$$
For general $N\ge 2$, by induction, \eqref{7} and \eqref{5}, it is readily seen that
\begin{align*}
|H^{\otimes N}(s,\bbx;t,\bby)|
&=\left|\int_s^t\int_{\mR^{2d}}H^{\otimes(N-1)}(s,\bbx;r,\bbz)b_1(r,\bbz)\cdot\nabla_{z_1}p_0(r,\bbz;t,\bby)\dif \bbz\dif r\right|\\
&\lesssim \Lambda_{\lambda_{N-2}}^{N-2}(t-s)|b_1(s,\bbx)|\int_s^t\!\!\int_{\mR^{2d}}g^{(1)}_{\l_{N-1}}(s,\bbx;r,\bbz)|b_1(r,\bbz)|g^{(1)}_{\lambda}(r,\bbz;t,\bby) \dif\bbz\dif r\\
&\le \Lambda_{\lambda_{N-2}}^{N-2}(t-s)|b_1(s,\bbx)|\int_s^t\!\!\int_{\mR^{2d}}g^{(1)}_{\l_{N-1}}(s,\bbx;r,\bbz)|b_1(r,\bbz)|g^{(1)}_{\lambda_{N-1}}(r,\bbz;t,\bby) \dif\bbz\dif r\\
&\lesssim \Lambda_{\lambda_{N-2}}^{N-2}(t-s) |b_1(s,\bbx)|\Lambda_{\lambda_{N-1}}(t-s) g^{(1)}_{\l_{N}}(s,\bbx;t,\bby)\\
&\leq\Lambda_{\lambda_{N-1}}^{N-1}(t-s)|b_1(s,\bbx)| g^{(1)}_{\l_{N}}(s,\bbx;t,\bby).
\end{align*}
Thus, we have for any $N\in\mN$,
\begin{align*}
|p_0\otimes H^{\otimes N}|(s,\bbx;t,\bby)
&\lesssim\Lambda_{\lambda_{N-1}}^{N-1}(t-s)\int_s^t\!\!\!\int_{\mR^{2d}}g^{(0)}_{\lambda}(s,\bbx;r,\bbz)|b_1(r,\bbz)|g^{(1)}_{\l_{N}}(r,\bbz;t,\bby)\dif\bbz\dif r\\
&\le\Lambda_{\lambda_{N-1}}^{N-1}(t-s)\int_s^t\!\!\!\int_{\mR^{2d}}g^{(0)}_{\l_{N}}(s,\bbx;r,\bbz)|b_1(r,\bbz)|g^{(1)}_{\l_{N}}(r,\bbz;t,\bby)\dif\bbz\dif r\\
&\lesssim \Lambda_{\lambda_{N-1}}^{N-1}(t-s)\Lambda_{\kappa\l_{N}}(t-s) g^{(0)}_{\kappa\l_{N}}(s,\bbx;t,\bby).
\end{align*}
This completes the proof.    
\end{proof}

To treat the remainder term $p\otimes H^{\otimes N}$, we use the same argument as in \cite{CMPZ23}.
We recall the following variational representation formula due to Bou\'e and Dupuis (see \cite{BD98} and \cite{Zh09}). 
\bt
Let $F$ be a bounded Wiener functional on the classical Wiener space $(\mC_T,\cB(\mC_T),\mu)$, where $\mu$ is the classical Wiener measure. 
Then we have
\begin{align}\label{28}
-\ln\mE^\mu \e^F=\inf_{h\in \cS_b}\mE^\mu \left(\frac{1}{2}\int_0^T|\dot{h}(s)|^2\dif s-F(\cdot+h)\right),
\end{align}
where $\cS_b$ denotes the set of all $\mR^d$-valued $\sF_t$-adapted and absolutely continuous process with 
\begin{align}\label{HH1}
\int_0^T|\dot{h}(s)|^2\dif s\leq N,\ \ N>0.
\end{align}
\et   

Let $h\in\cS_b$ be a control process satisfying \eqref{HH1}.
Consider the  following controlled SDE:
\begin{align}\label{10}
\dif \bbX^{h}_t=B\(t,\bbX^{h}_t\)\dif t+(\dot{h}(t),0)^*\dif t+\Sigma\(t,\bbX^{h}_t\)\dif W_t,\ \ \bbX^h_0=\bbx,
\end{align}
where $W$ is the Brownian motion on the classical Wiener space $(\mC_T,\cB(\mC_T),\mu)$. All the expectations $\mE$ below will be taken with respect to the classical Wiener measure $\mu$.
Since $B$ and $\Sigma$ are Lipschitz continuous, there is a unique solution $\bbX^h$ to the above controlled SDE. 
Still, we fix $T$ as in \eqref{13}.
We need the following conditional version of Krylov's estimate for controlled process $\bbX^h$.

\bl[Krylov type estimate]\label{Kry}
For any $h\in \cS_b$, there exist two constants $C=C(\Theta)>0$ and $\lambda\in(0,1)$ such that  for all $0\leq f\in\mL^\infty_TC^\infty_b(\mR^{2d})$ and $0\leq s<t\leq T$,
\begin{align}\label{9}
\mE\left(\int_s^tf(r,\bbX_r^h)\dif r\big|\sF_s\right)\lesssim_C K^{(1)}_\lambda(f;T)\left[\mE\left(\int_s^t|\dot{h}(r)|\dif r\big|\sF_s\right)+1\right].
\end{align}
\el
\begin{proof}
Let $u(s,\bbx):=\sI^T_f(s,\bbx)$.
By  \eqref{AG1}, \eqref{Sc1} and It\^o's formula, we have
\begin{align*}
u(t,\bbX_t^h)
&=u(s,\bbX_s^h)+\int_s^t(\p_r+\sL_r+\dot{h}(r)\cdot\nabla_{x_1})u(r,\bbX_r^h)\dif r+\int_s^t(\Sigma^*\cdot \nabla_{x_1}u)(r,\bbX_r^h)\dif W_r\\
&=u(s,\bbX_s^h)+\int_s^t[(b_1+\dot{h}(r))\cdot\nabla_{x_1}u(r,\bbX_r^h)-f(r,\bbX_r^h)]\dif r+\int_s^t(\Sigma^*\cdot \nabla_{x_1}u)(r,\bbX_r^h)\dif W_r.
\end{align*} 
By taking conditional expectation and \eqref{S6}, we have
\begin{align}
&\mE\left(\int_s^tf(r,\bbX_r^h)\dif r\big|\sF_s\right)\le 2\|u\|_{\mL^\infty_T}+\|\nabla_{x_1}u\|_{\mL^\infty_T}\mE\left(\int_s^t(|b_1(r,\bbX^h_r)|+|\dot{h}(r)|)\dif r\big|\sF_s\right)\no\\
&\qquad\le 2C_0K^{(0)}_{\lambda}(f; T)+C_1K^{(1)}_{\lambda}(f; T)\mE\left(\int_s^t(|b_1(r,\bbX^h_r)|+|\dot{h}(r)|)\dif r\big|\sF_s\right).\label{Kry1}
\end{align}
Taking $f=|b_1|$ in \eqref{Kry1} and by \eqref{13}, we get
$$
\mE\left(\int_s^t|b_1(r,\bbX_r^h)|\dif r\big|\sF_s\right)
\le 2C_0K^{(0)}_{\lambda}(|b_1|; T)+\frac12\mE\left(\int_s^t(|b_1(r,\bbX^h_r)|+|\dot{h}(r)|)\dif r\big|\sF_s\right),
$$
which implies that
\begin{align}
\mE\left(\int_s^t|b_1(r,\bbX_r^h)|\dif r\big|\sF_s\right)
&\le 4C_0K^{(0)}_{\lambda}(|b_1|; T)+ \mE\left(\int_s^t|\dot{h}(r)|\dif r\big|\sF_s\right)\no\\
&\stackrel{\eqref{AA1}, \eqref{13}}{\leq} C+ \mE\left(\int_s^t|\dot{h}(r)|\dif r\big|\sF_s\right).\label{AZ7}
\end{align} 
Substituting this into \eqref{Kry1} and by \eqref{AA1}, we obtain the desired estimate.
\end{proof}
Now we can show the following variational estimate that is crucial for the upper bound estimate.
\bl\label{Le54}
Let $\ell:\mR^{2d}\to(0,\infty)$ be a measurable function satisfying that for some $\kappa>1$,
$$
\kappa^{-1}\le \ell(\bbx) \le\kappa,\ \ \bbx\in\mR^{2d}.
$$
For the given $T>0$ from \eqref{13}, there is a constant $\lambda=\lambda(\Theta)>0$ such that for all $\bbx\in\mR^{2d}$ and $0\le s\le t\le T$,
\begin{align*}
\mE \ell(\bbX_{t,s}(\bbx))\le \sup_{\bbz\in\mR^{2d}}\e^{\ln \ell(\bbz)-\lambda|\bbz-\t_{t,s}(\bbx)|^2+\frac12}.
\end{align*}
\el
\begin{proof}
Without loss of generality, we assume $s=0$ and write 
$$
\bbX_t:=\bbX_{t,s}(x),\,\t_{t}=\t_{t,0}(\bbx).
$$ 
By the variational representation \eqref{28} with $F=\ln(\ell(\bbX_t))$, we have
\begin{align}\label{092}
-\ln\mE \ell(\bbX_t)=\inf_{h\in \cS_b}\mE\left(\frac{1}{2}\int_0^t|\dot{h}(s)|^2\dif s-\ln \ell(\bbX_t^h)\right),
\end{align}
where $\bbX^h$ solves SDE $\eqref{10}$ with starting point $\bbx$. Then by $\eqref{AZ7}$ we have 
\begin{align*}
\mE\left(\int_0^t|b_1(s,\bbX_s^h)|\dif s\right)^2
&=2\mE\left(\int_0^t|b_1(s,\bbX_s^h)|\mE\left(\int^t_s|b_1(r,\bbX_r^h)|\dif r\Big|\sF_s\right)\dif s\right)\\
&\leq  2\mE\left(\int_0^t|b_1(s,\bbX_s^h)|\mE\left(\int_s^t|\dot{h}(r)|\dif r\big|\sF_s+C\right)\dif s\right)\\
&\leq2 \mE\left(\int_0^t|b_1(s,\bbX_s^h)|\dif s\left(\int_0^t|\dot{h}(r)|\dif r+C\right)\right)\\
&\leq\frac12 \mE\left(\int_0^t|b_1(s,\bbX_s^h)|\dif s\right)^2+2\mE\left(\int_0^t|\dot{h}(r)|\dif r+C\right)^2,
\end{align*}
where the last step is due to Young's inequality. This implies that
\begin{align}\label{091}
\mE\left(\int_0^t|b_1(r,\bbX_r^h)|\dif r\right)^2\leq 4\mE\left(\int_0^t|\dot{h}(r)|\dif r+C\right)^2.
\end{align}
On the other hand, by \eqref{100}, \eqref{10} and  It\^o's formula, we have
\begin{align*}
|\bbX_t^{h}-\t_t|^2&=2\int_0^t\langle B(r,\bbX_r^h)+(\dot h(r),0)^*-B_0(r,\t_r),\bbX_r^{h}-\t_r\rangle\dif r\\
&\quad+2\int^t_0\<\bbX_r^{h}-\t_r, \Sigma(r, \bbX^h_r)\dif W_r\>+\int^t_0\tr(\Sigma\Sigma^*)(r,\bbX_r^h) \dif r.
\end{align*}
Noting that by \eqref{23},
\begin{align*}
|B(r,\bbx)-B_0(r,\bby)|\leq|b_0(r,\bbx)-b_0(r,\bby)|+|x_2-y_2|+|b_1(r,\bbx)|
\leq \kappa_1|\bbx-\bby|+|b_1(r,\bbx)|,
\end{align*}
by BDG's inequality and Young's inequality, we have for any $\eps>0$,
\begin{align*}
\mE\left(\sup_{s\in[0,t]}|\bbX_s^h-\t_s|^2\right)
&\leq \kappa_1\mE\left(\int_0^t|\bbX_r^h-\t_r|^2\dif r\right)+\mE\left(\int_0^t(|b_1(r,\bbX_r^h)|+|\dot h(r)|)\cdot|\bbX_r^h-\t_r|\dif r\right)\\
&\qquad+C\mE\left(\int_0^t|\bbX_r^h-\t_r|^2\dif r\right)^{1/2}+C t\\
&\leq\kappa_1\mE\left(\int_0^t|\bbX_r^h-\t_r|^2\dif r\right)+\eps\mE\left(\sup_{r\in[0,t]}|\bbX_r^h-\t_r|^2\right)\\
&\qquad+C_\eps\mE\left(\int_0^t(|b_1(r,\bbX_r^h)|+|\dot h(r)|)\dif r\right)^2+C_\eps t.
\end{align*}
Hence,
by Gronwall's inequality and \eqref{091},  we get
\begin{align*}
\mE\left(\sup_{s\in[0,t]}|\bbX_s^h-\t_s|^2\right)
\lesssim_{C_1} \mE\left(\int_0^t(|b_1(r,\bbX_r^h)|+|\dot h(r)|)\dif r\right)^2+t
\lesssim_{C_2}\mE\left(\int_0^t|\dot{h}(r)|^2\dif r\right)+1.
\end{align*}
Substituting this into \eqref{092}, we derive
\begin{align*}
-\ln\mE\ell(\bbX_t)&\ge \inf_{h\in \cS_b}\left(\frac{1}{2C_2}\mE|\bbX_t^h-\t_t| ^2-\frac12-\mE\ln \ell(\bbX_t^h)\right)\\
&=\inf_{h\in \cS_b}\mE\left(\frac{1}{2C_2}|\bbX_t^h-\t_t| ^2-\ln \ell(\bbX_t^h)\right)-\frac12\\
&\ge \inf_{\bbz\in\mR^{2d}}\left(\frac{1}{2C_2}|\bbz-\t_t| ^2-\ln \ell(\bbz)\right)-\frac12.
\end{align*}
From this we immediately obtain the desired estimate.
\end{proof}
\subsection{Proof of the upper bound}
Now we use the previous lemmas to show the following  Gaussian upper-bound estimate.
\bp[Control of the remainder]Suppose that
$|b_1|^\g\in\mK_1$ for some $\gamma>1$. Then
there exist constants $C_3,\lambda>0$ and $N_0\in\N$ such that for all $N\ge N_0$  and $(s,\bbx;t,\bby)\in\mD_T$,
$$
(p\otimes H^{\otimes N})(s,\bbx;t,\bby)\le C_3g_{\lambda}(s,\bbx;t,\bby).
$$
\ep
\begin{proof}
By a standard scaling technique, it suffices to consider the case $s=0,t=1$ (see \cite[Proposition 3.7]{CMPZ23} for more details about the scaling).
By $\eqref{5h}$ and H\"older's inqualiity we have
\begin{align*}
|(p\otimes H^{\otimes N})(0,\bbx;1,\bby)|&\lesssim\int_0^1 \Lambda^{N-1}_{\lambda_{N-1}}(1-r)\int_{\mR^{2d}}p(0,\bbx;r,\bbz)|b_1(r,\bbz)|g^{(1)}_{\l_N}(r,\bbz;1,\bby)\dif\bbz\dif r\\
&=\int_0^1\Lambda^{N-1}_{\lambda_{N-1}}(1-r)\mE\left(|b_1(r,\bbX_{r,0}(\bbx))|g^{(1)}_{\l_N}(r,\bbX_{r,0}(\bbx);1,\bby)\right)\dif r\\
&\le\left(\int_0^1\Lambda^{\gamma'(N-1)}_{\lambda_{N-1}}(1-r)\mE g_{\l_N}^{(1)}(r,\bbX_{r,0}(\bbx);1,\bby)^{\gamma'}\dif r\right)^\frac{1}{\gamma'}\\
&\quad\times \left(\mE\int_0^1 |b_1(r,\bbX_{r,0}(\bbx))|^\gamma\dif r \right)^\frac{1}{\gamma},
\end{align*}
where ${\frac1\g}+\frac{1}{\gamma'}=1$. By Lemma \ref{Kry} with $h=0$, we have
\begin{align*}
\mE\left( \int_0^1|b_1(r,\bbX_{r,0}(\bbx))|^\gamma\dif r\right)\lesssim K^{(1)}_\lambda(|b_1|^\gamma;1).
\end{align*}
On the other hand, for fixed $\eps,r\in(0,1)$ and $\bby\in\mR^{2d}$, let
$$
\ell_\eps(\bbx):=\eps+g^{(1)}_{\l_N}(r,\bbx;1,\bby)^{\gamma'}.
$$
Clearly, 
$$\eps\le\ell_\eps\le\eps+(1-r)^{-\gamma'(2d+1/2)}.$$
By Lemma \ref{Le54} with $\ell(\bbx)=\ell_\eps(\bbx)$ and \eqref{24}, \eqref{107}, we have
\begin{align*}
\mE \ell_\eps(\bbX_{r,0}(\bbx))&\le \sup_{\bbz\in\mR^{2d}}\e^{\ln \ell_\eps(\bbz)-\lambda|\bbz-\t_{r,0}(\bbx)|^2+\frac12}
=\sqrt{\e}\sup_{\bbz\in\mR^{2d}} \left(\ell_\eps(\bbz)\e^{-\lambda|\bbz-\t_{r,0}(\bbx)|^2}\right)\\
&\leq\sqrt{\e}\left(\eps+(1-r)^{-\gamma'(2d+1/2)}\sup_{\bbz\in\mR^{2d}} \e^{-\gamma'\l_N|\mT_{1-r}(\t_{1,r}(\bbz)-\bby)|^2-\lambda|\bbz-\t_{r,0}(\bbx)|^2}\right).
\end{align*}
Noting that by \eqref{16},
\begin{align*}
\sup_{\bbz\in\mR^{2d}} \e^{-\gamma'\l_N|\mT_{1-r}(\t_{1,r}(\bbz)-\bby)|^2-\lambda|\bbz-\t_{r,0}(\bbx)|^2}
&\le\sup_{\bbz\in\mR^{2d}} \e^{-\gamma'\l_N|\t_{1,r}(\bbz)-\bby|^2-\lambda|\bbz-\t_{r,0}(\bbx)|^2}\\
&\lesssim\sup_{\bbz\in\mR^{2d}} \e^{-\l'_N|\bbz-\theta_{r,1}(\bby)|^2-\lambda_N'|\bbz-\t_{r,0}(\bbx)|^2}\\
&\lesssim \e^{-\frac{\l'_N}2|\theta_{r,1}(\bby)-\t_{r,0}(\bbx)|^2}\lesssim \e^{-\l_N''|\t_{1}(\bbx)-\bby|^2},
\end{align*}
we further have
\begin{align*}
\mE g^{(1)}_{\l_N}(r,\bbX_{r,0}(\bbx);1,\bby)^{\gamma'}
&=\lim_{\eps\to 0}\mE \ell_\eps(\bbX_{r,0}(\bbx))\lesssim(1-r)^{-(2d+1/2)\gamma'}\e^{-\l_N''|\t_{1}(\bbx)-\bby|^2}.
\end{align*}
Hence,
\begin{align*}
|(p\otimes H^{\otimes N})(0,\bbx;1,\bby)|&\lesssim\left(\int_0^1\frac{\Lambda_{\lambda_{N-1}}(r)^{(N-1)\gamma'}}{r^{(2d+1/2)\gamma'}}\dif r\right)\e^{-\l_N''|\t_{1}(\bbx)-\bby|^2}.
\end{align*}
By \eqref{40} and Lemma \ref{KC}, there exists an  $N=N(\gamma,d)$ large enough so that
\begin{align*}
\int_0^1\frac{\Lambda_{\lambda_{N-1}}(r)^{(N-1)\gamma'}}{r^{(2d+1/2)\gamma'}}\dif r=\int_0^1\frac{K^{(1)}_{\lambda_{N-1}}(|b_1|;r)^{(N-1)\gamma'}}{r^{(2d+1/2)\gamma'}}\dif r<\infty.
\end{align*}
Thus, we have
$$
|(p\otimes H^{\otimes N})(0,\bbx;1,\bby)|\lesssim\e^{-\l_N''|\t_{1}(\bbx)-\bby|^2}.
$$
The proof is complete.
\end{proof}

\subsection{Proof of the lower bound}
The proof is almost the same as in \cite[Subsection 3.3]{CMPZ23}. For the readers' convenience, we provide the detailed proof here.
For fixed $(s,\bbx;t,\bby)\in\mD_T$,
consider the following deterministic control problem:
\begin{align}\label{ODE0}
\dot{\phi}_{r,s}=B_0(r,\phi_{r,s})+(\varphi_r,0)^*,\ r\in[s,t],\ \phi_{s,s}=\bbx,\ \phi_{t,s}=\bby,
\end{align}
where $\varphi:[s,t]\to\mR^{d}$ is a square integrable control function. Let $I(s,\bbx;t,\bby)$ be the associated energy functional 
\begin{align*}
I(s,\bbx;t,\bby)=\inf_{\varphi\in L^2([s,t];\mR^d)}\left\{\left(\int_s^t|\varphi_r|^2\dif r\right)^{1/2},\ \phi_{s,s}=\bbx,\ \phi_{t,s}=\bby\right\}.
\end{align*} 
We recall the following lemma proven in \cite[Proposition 2.5]{CMPZ23} about the energy functional which is crucial for proving the lower bound estimate.
\bl\label{19}
There exists constants $c_1=c_1(T,d,\kappa_1)\ge 1$ such that for all $(s,\bbx;t,\bby)\in\mD_T$,
\begin{align}\label{17}
c_1^{-1}(|\mT_{t-s}(\t_{t,s}(\bbx)-\bby)|-1)\le I(s,\bbx;t,\bby)\le c_1(|\mT_{t-s}(\t_{t,s}(\bbx)-\bby)|+1).
\end{align}
Moreover, one can find a control $\varphi\in L^2([s,t];\mR^d)$ and a solution $\phi_{r,s}$ to ODE \eqref{ODE0} such that 
\begin{align}\label{18}
\sup_{r\in[s,t]}|\varphi_r|\le c_2(|\mT_{t-s}(\t_{t,s}(\bbx)-\bby)|+1)/\sqrt{t-s},
\end{align}
where $c_2=c_2(T,d,\kappa_1)>0$.
\el
Now we can show the following lower bound estimate.
\bp
Suppose that
$|b_1|^\g\in\mK_1$ for some $\gamma>1$.  Then there exist constants $C,\l>0$ such that for all $(s,\bbx;t,\bby)\in\mD_T$,
$$
p(s,\bbx;t,\bby)\ge Cg_{\l}(s,\bbx;t,\bby).
$$
\ep
\begin{proof}
Let 
$
c_3=2c_1c_2+2c_0+2,
$
where $c_0,c_1$ and $c_2$ are from \eqref{16}, \eqref{17} and \eqref{18}. 
Since $\lim_{\delta\to 0}K^{(1)}_\lambda(|b_1|; \delta)=0$,  we choose $T_0>0$ such that
\begin{align}\label{211}
K^{(1)}_\lambda(|b_1|; T_0)\le\e^{-c^2_3/\lambda}/(2C_0C_3).
\end{align}
We first derive the lower bound on $\mD_{T_0}$. By \eqref{6}, \eqref{5}, \eqref{7} and the upper bound for $p$, we have
\begin{align*}
p(s,\bbx;t,\bby)&=p_0(s,\bbx;t,\bby)+(p\otimes H)(s,\bbx;t,\bby)\\
&\ge C_0^{-1}g^{(0)}_{\lambda^{-1}}(s,\bbx;t,\bby)-\int_s^t\!\!\!\int_{\mR^{2d}}p(s,\bbx;r,\bbz)|H(r,\bbz;t,\bby)|\dif\bbz\dif r\\
&\ge C_0^{-1}g^{(0)}_{\lambda^{-1}}(s,\bbx;t,\bby)-C_2\int_s^t\!\!\!\int_{\mR^{2d}}g^{(0)}_{\lambda}(s,\bbx;r,\bbz)|b_1(r,\bbz)|g^{(1)}_{\lambda}(r,\bbz;t,\bby)\dif\bbz\dif r\\
&\ge C^{-1}_0g^{(0)}_{\l^{-1}_0}(s,\bbx;t,\bby)-C_3K^{(1)}_\lambda(|b_1|; t-s)g^{(0)}_{\l_4}(s,\bbx;t,\bby)\\
&=(t-s)^{-2d}\Big(C^{-1}_0\e^{-|\mT_{t-s}(\t_{t,s}(\bbx)-\bby)|^2/\lambda}-C_3 K^{(1)}_\lambda(|b_1|; t-s)\e^{-\l_4|\mT_{t-s}(\t_{t,s}(\bbx)-\bby)|^2}\Big).
\end{align*}
If $(s,\bbx;t,\bby)\in\mD_T$ satisfies
\begin{align*}
|\mT_{t-s}(\t_{t,s}(\bbx)-\bby)|\le c_3,
\end{align*}
then by \eqref{211},
\begin{align}\label{21}
\begin{split}
p(s,\bbx;t,\bby)&\ge(t-s)^{-2d}\Big(C^{-1}_0\e^{-c^2_3/\lambda}-C_3K^{(1)}_\lambda(|b_1|; T_0)\Big)\\
&\ge\frac{\e^{-c_3^2/\lambda}}{2C_0}(t-s)^{-2d}=\hat{C}(t-s)^{-2d}.
\end{split}
\end{align}
Next we consider the case that
\begin{align*}
|\mT_{t-s}(\t_{t,s}(\bbx)-\bby)|\ge c_3.
\end{align*}
Let $M$ be the smallest integer such that 
\begin{align}\label{22}
M-1\le|\mT_{t-s}(\t_{t,s}(\bbx)-\bby)|^2\le M.
\end{align}
Define 
$$
\de:=\frac{t-s}{M}, \; t_j=s+j\de, \; j=0,1,\cdots, M
$$
Let $\phi_{t,s}(\bbx)$ be the optimal curve in Lemma \ref{19} and $\varphi$ is the corresponding control with
\begin{align*}
\sup_{r\in[s,t]}|\varphi_r|\le c_2(|\mT_{t-s}(\t_{t,s}(\bbx)-\bby)|+1)/\sqrt{t-s}\le c_2(\sqrt{M}+1)/\sqrt{t-s}.
\end{align*}
Define
$$
\xi_j:=\phi_{t_j,s}(\bbx),\;j=0,1,\cdots,M.
$$
For any $j=0,\cdots,M-1$, by Lemma \ref{19}, we have
\begin{align}\label{20}
|\mT_\delta(\t_{t_{j+1},t_j}(\xi_j))-\xi_{j+1}|
&\le c_1I(t_j,\xi_j;t_{j+1},\xi_{j+1})+ 1
\le c_1\left(\int_{t_j}^{t_{j+1}}|\varphi_r|^2\dif r \right)^{1/2}+ 1\nonumber\\
&\le c_1(t_{j+1}-t_j)^{1/2}\sup_{r\in[s,t]}|\varphi_r|+ 1\le  2c_1c_2+1.
\end{align}
Now set
$$
\S_0:=\{\xi_0\}=\{\bbx\},\;\; \S_M:=\{\xi_M\}=\{\bby\},
$$
and for $j=1,\cdots,M-1$,
$$
\S_j=\{\bbz\in\mR^{2d}:|\mT_\delta(\bbz-\xi_j)|\le 1\}.
$$
By \eqref{20} and  \eqref{16}, for any $j=0,\cdots,M-1$, we have that for $z_j\in\S_j$ and $z_{j+1}\in\S_{j+1}$,
\begin{align*}
|\mT_\delta(\t_{t_{j+1},t_j}(\bbz_j)-\bbz_{j+1})|
\le&|\mT_\delta(\t_{t_{j+1},t_j}(\xi_j)-\xi_{j+1})|+|\mT_\delta(\xi_{j+1}-\bbz_{j+1})|\\
&+|\mT_\delta(\t_{t_{j+1},t_j}(\bbz_j)-\t_{t_{j+1},t_j}(\xi_j))|\\
\le&|\mT_\delta(\t_{t_{j+1},t_j}(\xi_j)-\xi_{j+1})|+|\mT_\delta(\xi_{j+1}-\bbz_{j+1})|\\
&+c_0(|\mT_\delta(\xi_j-\bbz_j)|+1)\\
\le& 2c_1c_2+2c_0+2=c_3.
\end{align*}
Thus, by the Chapman-Kolmogorov equation and
\eqref{21}, we have
\begin{align*}
p(s,\bbx:t,\bby)
&=\int_{\mR^{2d}}\cdots\int_{\mR^{2d}}p(t_0,\bbx;t_1,\bbz_1)\cdots p(t_{M-1},\bbz_{M-1};t_M,\bby)\dif\bbz_1\cdots\dif\bbz_{M-1}\\
&\ge\int_{\S_1}\cdots\int_{\S_{M-1}}p(t_0,\bbx;t_1,\bbz_1)\cdots p(t_{M-1},\bbz_{M-1};t_M,\bby)\dif\bbz_1\cdots\dif\bbz_{M-1}\\
&\ge(\hat{C}\de^{-2d})^M\int_{\S_1}\cdots\int_{\S_{M-1}}\dif\bbz_1\cdots\dif\bbz_{M-1}=(\hat{C}\de^{-2d})^M(C\de^{2d})^{M-1},
\end{align*}
where $\hat{C}$ is given in \eqref{21}, and the last equality is due to $|\S_j|=C\de^{2d}$ for some $C$ depending on $d$. Recalling $\de=(t-s)/M$ and $M$ given in \eqref{22}, we finally have
\begin{align*}
p(s,\bbx;t,\bby)
&\ge\hat{C}^M C^{M-1}\de^{-2d}=(t-s)^{-2d}M^{2d}\exp\{M\log(\hat{C}C)\}/C\\
&\gtrsim (t-s)^{-2d}\exp\{-\lambda|\mT_{t-s}(\t_{t,s}(\bbx)-\bby)|^2\}.
\end{align*}
The lower bound is thus obtained for $t-s\le T_0$ and all $\bbx,\bby\in\mR^{2d}$. For general $0\le s\le t\le T$, it again follows from the Chapman-Kolmogorov equation. 
\end{proof}

{\bf Acknowledgement:} The authors would like to thanks Zimo Hao for quite useful conversation.

\end{document}